\documentclass[12pt]{article}
\usepackage{amsmath}
\usepackage{amssymb}
\usepackage{vatola}

\def\q{\quad}
\def\qq{\qquad}
\def\qtq#1{\q\t{#1}\q}
\def\mod#1{\ (\text{\rm mod}\ #1)}
\def\t{\text}
\def\f{\frac}
\def\e{\equiv}
\def\b{\binom}
\def\qp#1{q_p(#1)}
\def\sls#1#2{(\f{#1}{#2})}
 \def\ls#1#2{\big(\f{#1}{#2}\big)}
\def\Ls#1#2{\Big(\f{#1}{#2}\Big)}
\def\xp{\langle x\rangle_p}

\def\sumkp{\sum_{k=0}^{p-1}}
\let \pro=\proclaim
\let \endpro=\endproclaim

\begin{document}
 \centerline {\bf
Congruences for two types of Ap\'ery-like sequences}
\par\q\newline
\centerline{Zhi-Hong Sun}\newline \centerline{School of Mathematics
and Statistics}\centerline{Huaiyin Normal University}
\centerline{Huaian, Jiangsu 223300, P.R. China} \centerline{Email:
zhsun@hytc.edu.cn} \centerline{Homepage:
http://maths.hytc.edu.cn/szh1.htm}
 \abstract{In this paper we  present many results and
conjectures on congruences involving two types of Ap\'ery-like
sequences $\{G_n(x)\}$ and $\{V_n(x)\}$.
 \par\q
\newline MSC: Primary 11A07, Secondary
05A19,11B65,11B68,11E25,65Q30
 \newline Keywords: Ap\'ery-like numbers, congruence,
  recurrence relation, Bernoulli number, Euler number}
 \endabstract
\section*{1. Introduction}

\par \par For $s>1$ let $\zeta(s)=\sum_{n=1}^{\infty}\f 1{n^s}$.
In 1979, in order to prove that $\zeta(2)$ and
 $\zeta(3)$ are irrational,
 Ap\'ery [Ap] introduced the Ap\'ery numbers $\{A_n\}$ and $\{A'_n\}$
 given by
 $$ A_n= \sum_{k=0}^n\binom nk^2\binom{n+k}k^2\qtq{and}A'_n=\sum_{k=0}^n\b
 nk^2\b{n+k}k.$$
It is well known  (see [B2]) that
$$\align &(n+1)^3A_{n+1}=(2n+1)(17n(n+1)+5)A_n-n^3A_{n-1}\q (n\ge 1),
\\&(n+1)^2A'_{n+1}=(11n(n+1)+3)A'_n+n^2A'_{n-1}\q (n\ge 1).\endalign$$
\par
Let $\Bbb Z$ and $\Bbb Z^+$ be the set of integers and the set of
positive integers, respectively.  The Ap\'ery-like numbers $\{u_n\}$
of the first kind satisfy
$$u_0=1, \ u_1=\ b,\ (n+1)^3u_{n+1}
=(2n+1)(an(n+1)+b)u_n-cn^3u_{n-1}\q (n\ge 1),\tag 1.1$$ where
$a,b,c\in\Bbb Z$ and $c\not=0$. Let $[x]$ be the greatest integer
not exceeding $x$, and let
$$\align
&D_n=\sum_{k=0}^n\b nk^2\b{2k}k\b{2n-2k}{n-k},\q T_n=\sum_{k=0}^n\b
nk^2\b{2k}n^2, \\& b_n=\sum_{k=0}^{[n/3]}\b{2k}k\b{3k}k\b
n{3k}\b{n+k}k(-3)^{n-3k},
\\&V_n=\sum_{k=0}^n\b{2k}k^2\b{2n-2k}{n-k}^2=\sum_{k=0}^n
\b nk\b{n+k}k(-1)^k\b{2k}k^216^{n-k}.
\endalign$$
Then $\{A_n\}$, $\{D_n\}$, $\{b_n\}$, $\{T_n\}$ and $\{V_n\}$ are
Ap\'ery-like numbers of the first kind with
$(a,b,c)=(17,5,1),(10,4,64),(-7,-3,81),(12,4,16)$ and $(16,8,256)$,
respectively. The numbers $\{D_n\}$ are called Domb numbers, and
$\{b_n\}$ are called Almkvist-Zudilin numbers. For $\{A_n\}$,
$\{D_n\}$, $\{b_n\}$, $\{T_n\}$ and $\{V_n\}$ see A005259, A002895,
A125143, A290575 and A036917 in Sloane's database ``The On-Line
Encyclopedia of Integer Sequences". For the congruences involving
$T_n$ see the author's paper [S17], for the formulas and congruences
involving $V_n$, see [S20], [Su3], [Su5] and [W].
\par In 2009 Zagier [Z] studied the Ap\'ery-like numbers $\{u_n\}$
of the second kind  given by
$$u_0=1,\ u_1=b\qtq{and}(n+1)^2u_{n+1}=(an(n+1)+b)u_n-cn^2u_{n-1}\ (n\ge 1),\tag 1.2$$
where $a,b,c\in\Bbb Z$ and $c\not=0$. Let
 $$\align &f_n=\sum_{k=0}^n\b nk^3=\sum_{k=0}^n\b nk^2\b{2k}n,\q
 \\&S_n=\sum_{k=0}^{[n/2]}\b{2k}k^2\b n{2k}4^{n-2k}=\sum_{k=0}^n\b nk\b{2k}k\b{2n-2k}{n-k},
 \\&a_n=\sum_{k=0}^n\b nk^2\b{2k}k,
 \q Q_n=\sum_{k=0}^n\b nk(-8)^{n-k}f_k,
 \\&W_n=\sum_{k=0}^{[n/3]}\b{2k}k\b{3k}k\b n{3k}(-3)^{n-3k},
 \\&G_n=\sum_{k=0}^n\b{2k}k^2\b{2n-2k}{n-k}4^{n-k}
 =\sum_{k=0}^n\b nk(-1)^k\b{2k}k^216^{n-k}.\endalign$$
 According to [Z] and [AZ], $\{A'_n\},\ \{f_n\},\ \{S_n\},\ \{a_n\},\
 \{Q_n\}$, $\{W_n\}$ and $\{G_n\}$ are  Ap\'ery-like sequences of the second kind with
$(a,b,c)=(11,3,-1),(7,2,-8),(12,4,32),(10,3,9),(-17,$ $-6,72), (-9,$
$-3,27)$ and $(32,12,256)$, respectively. The  sequence
 $\{f_n\}$ is called Franel numbers.
 In [JS,S15,S16,S18,S20]  the author
 systematically investigated identities and congruences for sums
 involving $S_n$, $f_n$, $W_n$ and $G_n$.
  For $\{A'_n\},\ \{f_n\},\ \{S_n\},\ \{a_n\},\
 \{Q_n\}$, $\{W_n\}$ and $\{G_n\}$ see A005258, A000172,
 A081085, A002893, A093388,
A291898 and A143583 in Sloane's database ``The On-Line Encyclopedia
of Integer Sequences".
\par Ap\'ery-like numbers have fascinating properties and they are
concerned with Picard-Fuchs differential equation, modular forms,
hypergeometric series, elliptic curves, series for $\f 1{\pi}$,
supercongruences, binary quadratic forms, combinatorial identities,
Bernoulli numbers and Euler numbers. See for example [AO], [AZ],
[B1]-[B2], [CCS], [CZ], [CTYZ], [GMY], [SB] and [Z].
\par For a prime $p$ let $\Bbb Z_p$ be the set of
rational numbers whose denominator is not divisible by $p$, and let
$\ls ap$ be the Legendre symbol. For positive integers $a,b$ and
$n$, if $n=ax^2+by^2$ for some integers $x$ and $y$, we briefly
write that $n=ax^2+by^2$.
\par Inspired by [AZ], [Zu] and [Su7], in this paper we systematically
investigate the properties of two types of Ap\'ery-like sequences
$\{G_n(x)\}$ and $\{V_n(x)\}$, which are generalizations of $G_n$
and $V_n$.   Define
$$G_n(x)=\sum_{k=0}^n\b nk(-1)^k\b xk\b{-1-x}k\q(n=0,1,2,\ldots).\tag 1.3$$
Then $m^nG_n(x)$ is an Ap\'ery-like sequence of the second kind with
$a=2m$, $b=m(x^2+x+1)$ and $c=m^2$. Also,
$$G_n(x)=\sum_{k=0}^n\b xk^2(-1)^{n-k}\b{-1-x}{n-k}.$$
We note that Z.W. Sun [Su7] introduced $G_n(x)$ earlier, and showed
that for any prime $p>3$ and $x\in\Bbb Z_p$ with $x\not\e -\f 12\mod
p$,
$$\align &\sum_{n=0}^{p-1}G_n(x)^2\e (-1)^{\xp}p\f{1+2(x-\xp)/p}{1+2x}\mod
{p^2},\\&\sum_{n=0}^{p-1}(2n+1)G_n(x)^2\e 0\mod {p^2}, \endalign$$
where $\xp\in\{0,1,\ldots,p-1\}$ is given by $x\e \xp\mod p$.
  It
is well known (see [S5-S8]) that
$$\aligned &\b{-\f 12}k=\f{\b{2k}k}{(-4)^k},\q\b{-\f 13}k\b{-\f 23}k=
\f{\b{2k}k\b{3k}k}{27^k},
\\&\b{-\f 14}k\b{-\f 34}k=\f{\b{2k}k\b{4k}{2k}}{64^k},
\q \b{-\f 16}k\b{-\f
56}k=\f{\b{3k}k\b{6k}{3k}}{432^k}.\endaligned\tag 1.4$$ Thus,
$G_n=16^nG_n(-\f 12)$. Define
$$\align &G_n^{(3)}=27^nG_n\Big(-\f 13\Big)=\sum_{k=0}^n\b
nk(-1)^k\b{2k}k\b{3k}k27^{n-k},\tag 1.5
\\&G_n^{(4)}=64^nG_n\Big(-\f 14\Big)=\sum_{k=0}^n\b
nk(-1)^k\b{2k}k\b{4k}{2k}64^{n-k},\tag 1.6
\\&G_n^{(6)}=432^nG_n\Big(-\f 16\Big)=\sum_{k=0}^n\b
nk(-1)^k\b{3k}k\b{6k}{3k}432^{n-k}.\tag 1.7\endalign$$ Then
$G_n^{(3)}, G_n^{(4)}$ and $G_n^{(6)}$ are Ap\'ery-like sequences of
the second kind with $(a,b,c)=(54,21,729),(128,52,4096)$ and
$(864,372,186624)$, respectively. In Section 3, for any prime $p>3$
and $x\in\Bbb Z_p$ with $x\not\e 0,\pm 1,2\mod p$, we obtain
congruences for $\sum_{n=0}^{p-1}G_n(x)$ and
$\sum_{n=0}^{p-1}nG_n(x)$ modulo $p^3$. See Theorems 3.2 and 3.3. We
also get congruences for $G_p(x)$, $G_{p-1}(x)\mod {p^3}$ and
$G_{\f{p-1}2}(x)\mod {p^2}$. For instance, taking $x=-\f 16$ yields
$$\align &\sum_{n=0}^{p-1}\f{G_n^{(6)}}{432^n}\e
\cases p^2\mod {p^3} &\t{if $p\e 1\mod 4$,}
 \\-\f {67}5p^2\mod {p^3}
&\t{if $p\e 3\mod 4$,}\endcases
\\&\sum_{n=0}^{p-1}\f{nG_n^{(6)}}{432^n}\e \cases -\f 5{77}p^2\mod
{p^3} &\t{if $p\e 1\mod 4$,}
 \\\f {2567}{385}p^2\mod {p^3}
&\t{if $p\e 3\mod 4$ and $p\not=7,11$,}\endcases
\\&G_{p-1}^{(6)}\e (-1)^{\f{p-1}2}186624^{p-1}+\f{155}9p^2E_{p-3}\mod
{p^3},
\\&G_p^{(6)}\e 372+8640(-1)^{\f{p-1}2}p^2E_{p-3}\mod {p^3},
\\&G_{\f{p-1}2}^{(6)}\e \cases 432^{\f{p-1}2}\ls p3\cdot 4x^2-2p\mod {p^2}&\t{if $p\e
1\mod 4$  and so $p=x^2+4y^2$,}\\0\mod {p}&\t{if $p\e 3\mod 4$,}
\endcases
\endalign$$
where $\{E_n\}$ are Euler numbers given by
$$E_{2n-1}=0,\q E_0=1,\q E_{2n}=-\sum_{k=1}^{n} \b
{2n}{2k}E_{2n-2k}\q(n\ge 1).$$ We also present many conjectures on
congruences involving
 $G_n^{(3)},
G_n^{(4)}$ and $G_n^{(6)}$.
\par In Section 4, we introduce
$$V_n(x)=\sum_{k=0}^n\b nk\b{n+k}k(-1)^k\b
xk\b{-1-x}k\q(n=0,1,2,\ldots)\tag 1.8$$ and state that $m^nV_n(x)$
is an Ap\'ery-like sequence of the first kind with $a=m$,
$b=m(2x^2+2x+1)$ and $c=m^2$. We show that
$$V_n(x)=\sum_{k=0}^n\b xk^2\b{-1-x}{n-k}^2
=\sum_{k=0}^n\b nk\b{n+k}k(-1)^{n-k}G_k(x)\tag 1.9$$ and
$$\sum_{n=0}^{p-1}\f{V_n(x)}{m^n}\e
\Big(\sum_{k=0}^{p-1}\f{G_k(x)}{m^k}\Big)^2\mod p\tag 1.10$$ for any
prime $p>3$ and $m,x\in\Bbb Z_p$ with $m\not\e 0\mod p$, and
establish congruences for
$$\align &\sum_{n=0}^{p-1}V_n(x)\mod {p^4},
\q \sum_{n=0}^{p-1}(2n+1)V_n(x)\mod {p^5}, \q
\sum_{n=0}^{p-1}(2n+1)(-1)^nV_n(x)\mod {p^4},
\\&\sum_{n=0}^{p-1}(-1)^nV_n(x)\mod {p^2},\q V_p(x)\mod {p^3},
\q V_{p-1}(x)\mod {p^3}.\endalign$$ Clearly, $V_n=16^nV_n(-\f 12)$.
Define
$$\align &V_n^{(3)}=27^nV_n\Big(-\f 13\Big)=\sum_{k=0}^n\b
nk\b{n+k}k(-1)^k\b{2k}k\b{3k}k27^{n-k},\tag 1.11
\\&V_n^{(4)}=64^nV_n\Big(-\f 14\Big)=\sum_{k=0}^n\b
nk\b{n+k}k(-1)^k\b{2k}k\b{4k}{2k}64^{n-k},\tag 1.12
\\&V_n^{(6)}=432^nV_n\Big(-\f 16\Big)=\sum_{k=0}^n\b
nk\b{n+k}k(-1)^k\b{3k}k\b{6k}{3k}432^{n-k}.\tag 1.13\endalign$$ Then
$V_n^{(3)},V_n^{(4)}$ and $V_n^{(6)}$ are Ap\'ery-like sequences of
the first kind with $(a,b,c)=(27,15,729),$ $(64,40,4096)$ and
$(432,312,186624)$, respectively.  By [AZ] or [Zu],
$$V_n^{(4)}=\sum_{k=0}^n\b{2k}k^3\b{2n-2k}{n-k}16^{n-k}.\tag 1.14$$
In Section 4, we obtain many congruences involving $V_n^{(3)},
V_n^{(4)}$ and $V_n^{(6)}$. As typical examples, for any prime
$p>3$,
$$\align &V_p^{(6)}\e 312\mod {p^3},\q V_{p-1}^{(6)}
\e 186624^{p-1}\mod {p^3},
\\&\sum_{n=0}^{p-1}\f{V_n^{(3)}}{(-27)^n}\e
\cases 4x^2-2p\mod {p^2}&\t{if $p=x^2+3y^2\e 1\mod 3$,}
\\0\mod p&\t{if $p\e 2\mod 3$,}\endcases
\\&\sum_{n=0}^{p-1}\f{V_n^{(4)}}{64^n}\e
(-1)^{\f{p-1}2}p+13p^3E_{p-3}\mod {p^4},
\\&\sum_{n=0}^{p-1}(2n+1)\f{V_n^{(3)}}{27^n}\e
p^3+\f{21}2p^4\qp 3\mod {p^5},
\\&\sum_{n=0}^{p-1}(2n+1)\f{V_n^{(6)}}{(-432)^n}\e
(-1)^{\f{p-1}2}p+\f{155}9p^3E_{p-3}\mod {p^4},
\endalign$$
where $q_p(a)$ is the Fermat quotient given by
$q_p(a)=(a^{p-1}-1)/p$. We also make some conjectures on congruences
involving  $V_n^{(3)}, V_n^{(4)}$ and $V_n^{(6)}$.

\par In Section 5, we present many congruences and conjectures
involving $a_n$, $Q_n$ and other Ap\'ery-like sequences. Let $p$ be
an odd prime, $m\in\Bbb Z_p$ and $(m+2)(m-2)\not\e 0\mod p$. We show
that
$$\align &\sum_{k=0}^{p-1}\b{2k}k\f{a_k}{(m+2)^k}
\e\Ls{(m+2)(m-2)}p\sum_{k=0}^{p-1}\b{2k}k\f{f_k}{(m-2)^k}\mod p,
\\&\sum_{k=0}^{p-1}\b{2k}k\f{Q_k}{(-8(m+2))^k}
\e\Ls{(m+2)(m-2)}p\sum_{k=0}^{p-1}\b{2k}k\f{f_k}{(-8(m-2))^k}\mod p,
\\&\sum_{k=0}^{p-1}\b{2k}k\f{Q_k}{(-9(m+2))^k}
\e\Ls{(m+2)(m-2)}p\sum_{k=0}^{p-1}\b{2k}k\f{a_k}{(-9(m-2))^k}\mod p.
\endalign$$
From this we deduce many congruences modulo $p$ in terms of certain
binary quadratic forms.
\par   In addition to the above, throughout this paper
we use the following notation. For $n\in \Bbb Z^+$ define $H_n=1+\f
12+\cdots+\f 1n$ and $H_n^{(r)}=1+\f 1{2^r}+\cdots+\f 1{n^r}$ for
$r\in\Bbb Z^+$. For convenience, we also assume $H_0=H_0^{(r)}=0$.
The Bernoulli numbers $\{B_n\}$ and the sequence $\{U_n\}$ are
defined by
$$\align &B_0=1,\q\sum_{k=0}^{n-1}\b nkB_k=0\q(n\ge 2),
\\& U_{2n-1}=0,\q U_0=1,\q U_{2n}=-2\sum_{k=1}^{n}
\b {2n}{2k}U_{2n-2k}\q(n\ge 1).
\endalign$$
For congruences involving $B_n$, $E_n$ and $U_n$ see [S1,S2,S4].
\section*{2. Preliminaries}
For any nonnegative integer $m$ and real number $x$, it is clear
that
$$\align \sum_{k=0}^m\b xk(-1)^k&=\sum_{k=0}^m\b{x-1}k(-1)^k+
\sum_{k=1}^m\b{x-1}{k-1}(-1)^k
\\&=\sum_{r=0}^m\b{x-1}r(-1)^r-\sum_{r=0}^{m-1}\b{x-1}r(-1)^r
=\b{x-1}m(-1)^m.\endalign$$ That is,
 $$\sum_{k=0}^m\b xk(-1)^k=\sum_{k=0}^m\b{-x-1+k}k
=\b{m-x}m.\tag 2.1$$ From [G, (1.52)],
$$\sum_{n=k}^m\b nk=\b{m+1}{k+1}.\tag 2.2$$
 By [S17, (11)],
$$\sum_{n=k}^{p-1}n\b nk=\Big(\f{p^2+p}{k+2}-\f
p{k+1}\Big)\b{p-1}k.\tag 2.3$$ By (2.2) and (2.3),
$$\aligned &\sum_{n=k}^{p-1}(2n+1)\b{n+k}{2k}
=2\sum_{n=k}^{p-1}(n+k)\b{n+k}{2k}-(2k-1)\sum_{n=k}^{p-1}\b{n+k}{2k}
\\&=2\Big(\f{(p+k)(p+k+1)}{2k+2}-\f {p+k}{2k+1}\Big)\b{p+k-1}{2k}-\b
{p+k}{2k+1}=\f{p(p-k)}{k+1}\b{p+k}{2k}.\endaligned\tag 2.4$$ Using
induction one can prove the identity (see [Su2, (3.5)])
$$\sum_{n=k}^{p-1}(2n+1)(-1)^n\b{n+k}{2k}=(p-k)\b{p+k}{2k}.\tag 2.5$$
\par  For two sequences $\{u_n\}$ and $\{v_n\}$, we have the binomial
inversion formula:
$$v_n=\sum_{k=0}^n\b nku_k\ (n=0,1,2,\ldots)
\iff u_n=\sum_{k=0}^n\b nk(-1)^{n-k}v_k\ (n=0,1,2,\ldots).$$ By
[S14, Theorem 2.2],
$$\sum_{k=0}^n\b nk\b{n+k}k\Big(u_k-(-1)^{n-k}\sum_{r=0}^k\b
kru_r\Big)=0.\tag 2.6$$ Let $p$ be an odd prime. By [S8, Theorem
2.4], for $x\in\Bbb Z_p$,
$$\sum_{k=0}^{p-1}\b xk\b{-1-x}k\Big((-1)^{\xp}u_k-\sum_{r=0}^k
\b kr(-1)^ru_r\Big)\e 0\mod{p^2}.\tag 2.7$$ It is easy to prove (see
[S2, Lemma 2.9]) that for $k=1,2,\ldots,p-1$,
$$\b{p-1}k(-1)^k\e 1-pH_k+\f{p^2}2\big(H_k^2-H_k^{(2)}\big)\mod
{p^3}.\tag 2.8$$ In 1952 Ljunggren proved that for any prime $p>3$
and $m,n\in\Bbb Z^+$,
$$\b{mp}{np}\e \b mn\mod{p^3}.\tag 2.9$$

\par The Bernoulli polynomials $\{B_n(x)\}$ and Euler polynomials
$\{E_n(x)\}$ are given by
$$B_n(x)=\sum_{k=0}^n\b nkB_kx^{n-k}\qtq{and}
E_n(x)=\f 1{2^n}\sum_{k=0}^n\b nk(2x-1)^{n-k}E_k.\tag 2.10$$ Suppose
that $p$ is an odd prime. By [S1, Lemma 3.2],
 $$H_{\xp}^{(2)}=\sum_{r=1}^{\xp}\f 1{r^2}\e \sum_{r=1}^{\xp}r^{p-3}
 \e -\f{B_{p-2}(-x)}{p-2}\e \f 12B_{p-2}(-x)\mod p.\tag 2.11$$
By [S13, Lemma 2.2 (with $k=2$)],
$$\sum_{k=1}^{\xp}\f{(-1)^k}{k^2}\e\f 12(-1)^{\xp}E_{p-3}(-x)\mod
p.\tag 2.12$$ By [S13, Theorem 2.1], for $x\in\Bbb Z_p$ and
$x'=(x-\xp)/p$,
$$\sum_{k=0}^{p-1}\b xk\b{-1-x}k\e (-1)^{\xp}+p^2x'(x'+1)E_{p-3}(-x)
\mod {p^3}.\tag 2.13$$
 By [S13, pp.3300-3301], for $p>3$,
$$\align &E_{p-3}\Ls 12=\f 1{2^{p-3}}E_{p-3}\e 4E_{p-3}\mod p,
\q E_{p-3}\Ls 16\e 20E_{p-3}\mod p,\tag 2.14
\\&E_{p-3}\Ls 13\e 9U_{p-3}\mod p,\q E_{p-3}\Ls 14\e 16s_{p-3}\mod
p,\tag 2.15
\endalign$$
where $\{s_n\}$ is given by
$$s_0=1\q \t{and}\q s_n=1-\sum_{k=0}^{n-1} \b
nk2^{2n-1-2k}s_k\q (n\ge 1).\tag 2.16$$ From [S10, Lemma 2.3 and
Theorem 2.1], for $x\in\Bbb Z_p$ and $x'=(x-\xp)/p$,
 $$  \sum_{k=1}^{p-1}\f{\b
 xk\b{-1-x}k}k\e -2H_{\xp}+2px'H_{\xp}^{(2)}
 \e -2\f{B_{p^2(p-1)}(-x)-B_{p^2(p-1)}}{p^2(p-1)}
 \mod {p^2}.\tag 2.17$$
 By [T1],
$$\sum_{k=1}^{p-1}\f{\b xk\b{-1-x}k}{k^2}\e -\f 12
\Big(\sum_{k=1}^{p-1}\f{\b xk\b{-1-x}k}k\Big)^2\mod p.\tag 2.18$$ By
[S10, (2.5)] and [S1, Theorem 5.2],
$$\f{B_{p^2(p-1)}(\f 12)-B_{p^2(p-1)}}{p^2(p-1)}
\e H_{\f{p-1}2}+\f p2H_{\f{p-1}2}^{(2)}\e -2\qp 2+p\qp 2^2\mod
{p^2}.\tag 2.19$$ By [S2, p.287], for $p>3$,
$$\align &\f{B_{p^2(p-1)}(\f 13)-B_{p^2(p-1)}}{p^2(p-1)}
\e-\f 32\qp 3+\f 34p\qp 3^2\mod {p^2},\tag 2.20
\\&\f{B_{p^2(p-1)}(\f 14)-B_{p^2(p-1)}}{p^2(p-1)}
\e -3\qp 2+\f 32p\qp 2^2\mod {p^2},\tag 2.21
\\&\f{B_{p^2(p-1)}(\f 16)-B_{p^2(p-1)}}{p^2(p-1)}
\e -2\qp 2-\f 32\qp 3+p\Big(\qp 2^2+\f 34\qp 3^2\Big)\mod {p^2}.\tag
2.22\endalign$$

\pro{Lemma 2.1} Let $p$ be an odd prime and $m\in\Bbb Z_p$ with
$m(m-1)\not\e 0\mod p$. Suppose that $u_0,u_1,\ldots,u_{p-1}\in\Bbb
Z_p$ and $v_n=\sum_{k=0}^n\b nku_k$ $(n\ge 0)$. Then
$$\sum_{k=0}^{p-1}\f{v_k}{m^k}
\e \sum_{k=0}^{p-1}\f{u_k}{(m-1)^k}\mod p.$$ Also,
$$\sum_{k=0}^{p-1}v_k\e \sum_{k=0}^{p-1}\f{(-1)^kp}{k+1}u_k
-p^2\sum_{k=0}^{p-2}\f{(-1)^kH_k}{k+1}u_k \mod {p^3}.$$
\endpro
Proof. It is clear that
$$\align \sum_{k=0}^{p-1}\f{v_k}{m^k}&
=\sum_{k=0}^{p-1}\f 1{m^k}\sum_{s=0}^k\b ksu_s =\sum_{s=0}^{p-1}
\sum_{k=s}^{p-1}\f 1{m^k}\b{-1-s}{k-s}(-1)^{k-s}u_s
\\&=\sum_{s=0}^{p-1}\f{u_s}{m^s}\sum_{k=s}^{p-1}\b{-1-s}{k-s}\f
1{(-m)^{k-s}}=\sum_{s=0}^{p-1}\f {u_s}{m^s} \sum_{r=0}^{p-1-s}
\b{-1-s}r\Big(-\f 1m\Big)^r
\\&\e \sum_{s=0}^{p-1}\f {u_s}{m^s} \sum_{r=0}^{p-1-s}\b{p-1-s}r
\Big(-\f 1m\Big)^r=\sum_{s=0}^{p-1}\f {u_s}{m^s}\Big(1-\f
1m\Big)^{p-1-s}
\\&\e \sum_{s=0}^{p-1}\f{u_s}{(m-1)^s}\mod p.\endalign$$
For $m=1$, from the above and (2.1) we see that
$$\align \sum_{k=0}^{p-1}v_k&=\sum_{s=0}^{p-1}u_s
\sum_{r=0}^{p-1-s} \b{-1-s}r(-1)^r =\sum_{s=0}^{p-1}\b p{p-1-s}u_s
\\&=\sum_{s=0}^{p-1}\f p{s+1}\b{p-1}su_s
\e u_{p-1}+\sum_{s=0}^{p-2}\f p{s+1}(-1)^s(1-pH_s)u_s \mod {p^3}.
\endalign$$
This yields the remaining part.

\par\q
\pro{Lemma 2.2} Let $p$ be an odd prime and $m\in\Bbb Z_p$ with
$(m+2)(m-2)\not\e 0\mod p$. Suppose that
$u_0,u_1,\ldots,u_{\f{p-1}2}\in\Bbb Z_p$ and $v_n=\sum_{k=0}^n\b
nku_k$ $(n\ge 0)$. Then
$$\sum_{k=0}^{(p-1)/2}\b{2k}k\f{v_k}{(m+2)^k}
\e\Ls{(m+2)(m-2)}p\sum_{k=0}^{(p-1)/2}\b{2k}k\f{u_k}{(m-2)^k}\mod
p\tag 2.23$$ and
$$\sum_{k=0}^{(p-1)/2}\b{2k}k\f{v_k}{4^k}
\e (-1)^{\f{p-1}2}\Big(2^{p-1}u_{\f{p-1}2}+\f
p2\sum_{s=0}^{(p-3)/2}\f{\b{2s}s}{(-4)^s(2s+1)}u_s\Big)\mod{p^2}.$$
\endpro
Proof. Note that $\b{-\f 12}k=\b{2k}k4^{-k}$ and $\b xk\b ks=\b
xs\b{x-s}{k-s}$. For $m\not\e -2\mod p$,
$$\align
&\sum_{k=0}^{(p-1)/2}\b{2k}k\f{v_k}{(m+2)^k}
\\&=\sum_{k=0}^{(p-1)/2}\b{-\f 12}k\Ls{-4}{m+2}^k \sum_{s=0}^k\b ksu_s
=\sum_{s=0}^{(p-1)/2}\sum_{k=s}^{(p-1)/2}\b{-\f 12}s \b{-\f
12-s}{k-s}u_s\Ls{-4}{m+2}^k
\\&=\sum_{s=0}^{(p-1)/2}\b{-\f 12}su_s\Ls{-4}{m+2}^s\sum_{k=s}^{(p-1)/2}
\b{-\f 12-s}{k-s}\Ls{-4}{m+2}^{k-s}
\\&=\sum_{s=0}^{(p-1)/2}\b{-\f 12}su_s\Ls{-4}{m+2}^s\sum_{r=0}^{(p-1)/2-s}
\b{-\f 12-s}r\Ls{-4}{m+2}^r.\endalign$$ Hence, for $m\not\e \pm
2\mod p$,
$$\align
&\sum_{k=0}^{(p-1)/2}\b{2k}k\f{v_k}{(m+2)^k}\\ &\e
\sum_{s=0}^{(p-1)/2}\b{-\f
12}su_s\Ls{-4}{m+2}^s\sum_{r=0}^{(p-1)/2-s} \b{\f
{p-1}2-s}r\Ls{-4}{m+2}^r
\\&=\sum_{s=0}^{(p-1)/2}\b{2s}s\f{u_s}{(m+2)^s}\Big(1-\f
4{m+2}\Big)^{\f{p-1}2-s}
=\sum_{s=0}^{(p-1)/2}\b{2s}s\f{u_s}{(m-2)^s}\Ls{m-2}{m+2}^{\f{p-1}2}
\\&\e \Ls{(m-2)(m+2)}p\sum_{s=0}^{(p-1)/2}\b{2s}s\f{u_s}{(m-2)^s}\mod
p.\endalign$$
 On the other hand, from (2.1), Morley's congruence
  $\b{p-1}{\f{p-1}2}
\e 4^{p-1}\mod {p^2}$ and the above we deduce that
$$\align &\sum_{k=0}^{(p-1)/2}\b{2k}k\f{v_k}{4^k}
\\&=\sum_{s=0}^{(p-1)/2}\b{-\f 12}su_s(-1)^s\sum_{r=0}^{(p-1)/2-s}
\b{-\f 12-s}r(-1)^r =\sum_{s=0}^{(p-1)/2}\f{\b{2s}s}{4^s}u_s\b{\f
p2}{\f{p-1}2-s} \\&=\f{\b{p-1}{(p-1)/2}}{4^{(p-1)/2}}u_{\f{p-1}2} +
\sum_{s=0}^{(p-1)/2}\f{\b{2s}s}{4^s}u_s\f{\f p2(\f p2-1)\cdots(\f
p2-(\f{p-3}2-s))}{(\f{p-1}2-s)!}
\\&\e (-1)^{\f{p-1}2}2^{p-1}u_{\f{p-1}2}+(-1)^{\f{p-1}2}p\sum_{s=0}^{(p-3)/2}
\f{\b{2s}su_s}{(-4)^s(2s+1)}\mod{p^2}.\endalign
$$
This proves the lemma.
\par{\bf Remark 2.1} (2.23) has been given by Z.W. Sun in [Su6, Theorem 2.2].

\pro{Lemma 2.3 ([S15, Lemma 2.4])} Let $p$ be an odd prime, $u,
c_0,c_1,\ldots,
 c_{p-1}\in\Bbb Z_p$ and
 $u\not\e  1\mod p$. Then
 $$\sum_{k=0}^{p-1}\b{2k}k\Ls u{(1-u)^2}^kc_k
 \e\sum_{n=0}^{p-1}u^n\sum_{k=0}^n\b nk\b {n+k}kc_k\mod p.$$
 \endpro
\pro{Lemma 2.4 ([S16, Lemma 2.2])} Let $p>3$ be a prime and
$c_0,c_1,\ldots,c_{p-1}\in\Bbb Z_p$. Then
$$\sum_{n=0}^{p-1}\sum_{k=0}^n\b nk\b{n+k}kc_k
\e \sum_{k=0}^{p-1}\f{p}{2k+1}(-1)^kc_k\mod{p^3}.$$
\endpro

\section*{3. Congruences for  $G_n(x)$}
\par Recall that
 $$G_n(x)=\sum_{k=0}^n\b nk(-1)^k\b xk\b{-1-x}k\q(n=0,1,2,\ldots).$$
Applying the binomial inversion formula,
$$\sum_{k=0}^n\b nk(-1)^kG_k(x)=\b xn\b{-1-x}n.$$
  Using sumtools in Maple we find that
$$\aligned&(n+1)^2m^{n+1}G_{n+1}(x)
\\&=(2mn(n+1)+m(x^2+x+1))m^nG_n(x) -m^2n^2m^{n-1}G_{n-1}(x)\ (n\ge
1)\endaligned\tag 3.1$$
 Thus, $m^nG_n(x)$ is an Ap\'ery-like
sequence of the second kind with $a=2m$, $b=m(x^2+x+1)$ and $c=m^2$.

\pro{Theorem 3.1} For $n=0,1,2,\ldots$ we have
$$G_n(x)=\sum_{k=0}^n\b xk^2(-1)^{n-k}\b{-1-x}{n-k}$$
\endpro
Proof. Set $G'_n(x)=\sum_{k=0}^n\b xk^2(-1)^{n-k}\b{-1-x}{n-k}$.
Then $G'_0(x)=1=G_0(x)$ and $G'_1(x)=x^2+x+1=G_1(x)$. Using sumtools
in Maple we find that
$$(n+1)^2G'_{n+1}(x)=(2n(n+1)+x^2+x+1)G'_n(x)-n^2G'_{n-1}(x)\q(n\ge
1).$$ Thus, $G'_n(x)=G_n(x)$ for $n=0,1,2,\ldots$ by (3.1).
\par\q
\par Let $G_n^{(3)}, G_n^{(4)}$ and $G_n^{(6)}$ be given by
(1.5)-(1.7). Then clearly they are  Ap\'ery-like integral sequences.

 \pro{Theorem 3.2} Let $p$ be an odd prime, $x\in\Bbb Z_p$, $x\not\e
0,-1\mod p$ and $x'=(x-\xp)/p$. Then
$$\sum_{n=0}^{p-1}G_n(x)\e
 p^2\f{x'(x'+1)+1-(-1)^{\xp}} {x(x+1)}\mod {p^3}.$$
Hence, for $p>3$,
 $$\align &\sum_{n=0}^{p-1}\f{G_n^{(3)}}{27^n}\e
\cases p^2\mod {p^3} &\t{if $p\e 1\mod 3$,}
\\-8p^2\mod {p^3} &\t{if $p\e 2\mod 3$,}
\endcases
\\&\sum_{n=0}^{p-1}\f{G_n^{(4)}}{64^n}\e
\cases p^2\mod {p^3} &\t{if $p\e 1,3\mod 8$,}
\\-\f{29}3p^2\mod {p^3} &\t{if $p\e 5,7\mod 8$,}\endcases
\\&\sum_{n=0}^{p-1}\f{G_n^{(6)}}{432^n}\e
\cases p^2\mod {p^3} &\t{if $p\e 1\mod 4$,}
 \\-\f {67}5p^2\mod {p^3}
&\t{if $p\e 3\mod 4$.}\endcases\endalign$$
\endpro
Proof. Using (2.2),
$$\align \sum_{n=0}^{p-1}G_n(x)&=\sum_{n=0}^{p-1}\sum_{k=0}^n
\b nk(-1)^k\b xk\b{-1-x}k
\\&=\sum_{k=0}^{p-1}\b xk\b{-1-x}k(-1)^k\sum_{n=k}^{p-1}\b nk
\\&=\sum_{k=0}^{p-1}\b xk\b{-1-x}k(-1)^k\b p{k+1}
\\&=\sum_{k=0}^{p-1}\b xk\b{-1-x}k(-1)^k\f p{k+1}\b{p-1}k.
\endalign$$
Since $x\not\e p,p-1\mod p$, we have $\b x{p-1}\b{-1-x}{p-1}\e
\b{\xp}{p-1}\b{p-1-\xp}{p-1}\e 0\mod {p^2}$. Note that
$\b{p-1}k(-1)^k\e 1-pH_k\mod {p^2}$ for $k=0,1,\ldots,p-1$. We then
get
$$\sum_{k=0}^{p-1}G_k(x)
\e \sum_{k=0}^{p-1}\b xk\b {-1-x}k\f p{k+1} -p^2\sum_{k=0}^{p-1} \b
xk\b{-1-x}k\f{H_k}{k+1} \mod {p^3}.$$ By [S8, Lemma 2.4] and [S12,
Lemma 2.2],
$$\aligned\sum_{k=0}^{p-1}\b xk\b {-1-x}k\f p{k+1}&=
\b{x-1}{p-1}\b{-2-x}{p-1}=\f{p+x}{x+1}\b{x-1}{p-1}\b{-x-1}{p-1}
\\&\e\f{p+x}{x+1}\cdot\f {p^2x'(x'+1)}{x^2} \e p^2\f{x'(x'+1)}{x(x+1)}
\mod {p^3}.\endaligned\tag 3.2$$ Using the symbolic summation
package Sigma in Mathematica, Liu and Ni [LN] found the identity
$$\sum_{k=0}^n\b
nk\b{n+k}k(-1)^k\f{H_k}{k+1}=\f{(-1)^n-1}{n(n+1)}.\tag 3.3$$ Hence
$$\align &\sum_{k=0}^{p-1} \b
xk\b{-1-x}k\f{H_k}{k+1}\\&\e\sum_{k=0}^{p-2} \b
xk\b{x+k}k(-1)^k\f{H_k}{k+1}\e \sum_{k=0}^{p-2} \b
{\xp}k\b{\xp+k}k(-1)^k\f{H_k}{k+1} \\&=\sum_{k=0}^{\xp} \b
{\xp}k\b{\xp+k}k(-1)^k\f{H_k}{k+1}=\f{(-1)^{\xp}-1}{\xp(\xp+1)} \e
\f{(-1)^{\xp}-1}{x(x+1)}\mod p.\endalign$$ Therefore,
$$\sum_{k=0}^{p-1}G_k(x)\e
p^2\f{x'(x'+1)}{x(x+1)}-p^2\f{(-1)^{\xp}-1}{x(x+1)}\mod{p^3}.$$ Now
assume that $p>3$. For $m=3,4,6$ we see that $(-1)^{\langle-\f
1m\rangle_p}=(-1)^{[\f pm]}$. Recall that
$$ \f{G_k^{(3)}}{27^k}=G_k\Big(-\f 13\Big), \q\f{G_k^{(4)}}{64^k}=
G_k\Big(-\f 14\Big), \q  \f{G_k^{(6)}}{432^k}=G_k\Big(-\f 16\Big).$$
Taking $x=-\f 13,-\f 14,-\f 16$ in the above congruence for
$\sum_{k=0}^{p-1}G_k(x)\mod{p^3}$ yields the remaining part.
\par\q
\par{\bf Remark 3.1} Let $p$ be an odd prime. In [S20], the author conjectued
 that
$$\sum_{n=0}^{p-1}\f{G_n}{16^n}\e
  (4(-1)^{\f{p-1}2}-3)p^2\mod {p^3},$$
which was solved by Liu and Ni in [LN]. This congruence is an easy
consequence of Theorem 3.2 (with $x=-\f 12$).
\par\q
\pro{Theorem 3.3} Let $p$ be an odd prime, $x\in\Bbb Z_p$,
$(x-1)x(x+1)(x+2)\not\e 0\mod p$ and $x'=(x-\xp)/p$. Then
$$\sum_{n=0}^{p-1}nG_n(x)\e
\cases p^2\f{x'(x'+1)(1-x(x+1))-x(x+1)}{(x-1)x(x+1)(x+2)} \mod {p^3}
&\t{if $2\mid \xp$,}
\\p^2\f{x'(x'+1)(1-x(x+1))-(x-1)(x+2)}{(x-1)x(x+1)(x+2)} \mod {p^3}
&\t{if $2\nmid \xp$.}
\endcases$$
Hence, for $p>3$,
 $$\align &\sum_{n=0}^{p-1}\f{nG_n}{16^n}\e
 \cases -\f 19p^2\mod {p^3}
&\t{if $p\e 1\mod 4$,}
\\\f{31}9p^2\mod {p^3}&\t{if $p\e 3\mod
4$,}\endcases
\\&\sum_{n=0}^{p-1}\f{nG_n^{(3)}}{27^n}\e
\cases -\f 1{10}p^2\mod {p^3} &\t{if $p\e 1\mod 3$,}
\\\f{79}{20}p^2\mod {p^3} &\t{if $p\e 2\mod 3$ and $p\not=5$,}
\endcases
\\&\sum_{n=0}^{p-1}\f{nG_n^{(4)}}{64^n}\e
\cases -\f 3{35}p^2\mod {p^3} &\t{if $p\e 1,3\mod 8$,}
\\\f{503}{105}p^2\mod {p^3} &\t{if $p\e 5,7\mod 8$ and $p\not=5,7$,}\endcases
\\&\sum_{n=0}^{p-1}\f{nG_n^{(6)}}{432^n}\e
\cases -\f 5{77}p^2\mod {p^3} &\t{if $p\e 1\mod 4$,}
 \\\f {2567}{385}p^2\mod {p^3}
&\t{if $p\e 3\mod 4$ and $p\not=7,11$.}\endcases\endalign$$
\endpro
Proof. By the definition of $G_n(x)$ and (2.3),
$$\align \sum_{n=0}^{p-1}nG_n(x)&=\sum_{n=0}^{p-1}n\sum_{k=0}^n
\b nk(-1)^k\b xk\b{-1-x}k
\\&=\sum_{k=0}^{p-1}\b xk\b{-1-x}k(-1)^k\sum_{n=k}^{p-1}n\b nk
\\&=\sum_{k=0}^{p-1}\b xk\b{-1-x}k(-1)^k\Big(\f{p^2+p}{k+2}-\f
p{k+1}\Big)\b{p-1}k .
\endalign$$
For $k=p-2$ or $p-1$, we have $\b xk\e \b{\xp}k=0\mod p$. Thus, $\b
xk\b{-1-x}k\f 1{k+1},\b xk\b{-1-x}k\f 1{k+2}\in\Bbb Z_p$ for
$k=0,1,\ldots,p-1$. Since $\b{p-1}k(-1)^k\e 1-pH_k\mod {p^2}$, we
then have
$$\align \sum_{n=0}^{p-1}nG_n(x)&
\e \sum_{k=0}^{p-1}\b xk\b{-1-x}k\Big(\f{p^2+p}{k+2}-\f
p{k+1}\Big)(1-pH_k)
\\&\e (p^2+p)\sum_{k=0}^{p-1}\b xk\b{-1-x}k\f 1{k+2}
-p\sum_{k=0}^{p-1}\b xk\b{-1-x}k\f 1{k+1}
\\&\qq+p^2\sum_{k=0}^{p-1}\b xk\b{-1-x}k\f {H_k}{k+1}
-p^2\sum_{k=0}^{p-1}\b xk\b{-1-x}k\f {H_k}{k+2}\mod {p^3}.
\endalign$$
By [S19, Theorem 2.1],
$$\sum_{k=0}^{p-1}\b xk\b{-1-x}k\f 1{k+2}
\e\f{-px'(x'+1)}{4\b{\xp}2\b{p-1-\xp}2} \e
-\f{px'(x'+1)}{(x-1)x(x+1)(x+2)}\mod {p^2}.\tag 3.4$$ By the proof
of Theorem 3.2,
$$\align &p\sum_{k=0}^{p-1}\b xk\b {-1-x}k\f 1{k+1}
\e p^2\f{x'(x'+1)}{x(x+1)}\mod {p^3},
\\&\sum_{k=0}^{p-1}\b xk\b {-1-x}k\f{H_k}{k+1}
\e\f{(-1)^{\xp}-1}{x(x+1)}\mod p.\endalign$$ Using summation package
Sigma in Mathematica one can find and prove that for $n\ge 2$,
$$\sum_{k=1}^n\b nk\b{n+k}k(-1)^k\f{H_k}{k+2}
=\cases \f 1{(n-1)(n+2)}&\t{if $n$ is even,}\\-\f 1{n(n+1)}&\t{if
$n$ is odd.} \endcases\tag 3.5$$ Hence,
$$\align & \sum_{k=0}^{p-1}\b xk\b{-1-x}k\f {H_k}{k+2}
\\&=\sum_{k=0}^{p-1}\b xk\b{x+k}k(-1)^k\f {H_k}{k+2}
\e \sum_{k=0}^{\xp}\b{\xp}k\b{\xp+k}k(-1)^k\f {H_k}{k+2}
\\&=\cases \f 1{(\xp-1)(\xp+2)}\e \f 1{(x-1)(x+2)}\mod p
&\t{if $\xp$ is even,}\\-\f 1{\xp(\xp+1)}\e -\f 1{x(x+1)}\mod
p&\t{if $\xp$ is odd.}
\endcases\endalign$$
If $\xp$ is even, combining the above gives
$$\align \sum_{n=0}^{p-1}nG_n(x)
&\e -(p^2+p)\f{px'(x'+1)}{(x-1)x(x+1)(x+2)}-p^2\f{x'(x'+1)}{x(x+1)}
-p^2\f 1{(x-1)(x+2)}
\\&\e p^2\f{x'(x'+1)(1-x(x+1))+x(x+1)}{(x-1)x(x+1)(x+2)} \mod {p^3}.
\endalign$$
If $\xp$ is odd, from the above we deduce that
$$\align \sum_{n=0}^{p-1}nG_n(x)
&\e -(p^2+p)\f{px'(x'+1)}{(x-1)x(x+1)(x+2)}-p^2\f{x'(x'+1)}{x(x+1)}
-\f {2p^2}{x(x+1)}+\f{p^2}{x(x+1)}
\\&\e p^2\f{x'(x'+1)(1-x(x+1))-(x-1)(x+2)}{(x-1)x(x+1)(x+2)} \mod {p^3}.
\endalign$$
Taking $x=-\f 12,-\f 13,-\f 14,-\f 16$ in the above congruences
yields the remaining results.

\pro{Theorem 3.4} Let $p$ be an odd prime, $x\in\Bbb Z_p$ and
$x'=(x-\xp)/p$. Then
$$\align G_{p-1}(x)&\e
(-1)^{\xp}\Big(1-2p\f{B_{p^2(p-1)}(-x)-B_{p^2(p-1)}}{p^2(p-1)}+2p^2
\Ls {B_{p^2(p-1)}(-x)-B_{p^2(p-1)}}{p^2(p-1)}^2\Big)
\\&\qq+p^2(x'(x'+1)+1)E_{p-3}(-x)\mod {p^3}.\endalign$$
 Hence, for $p>3$,
 $$\align &G_{p-1}^{(3)}\e (-1)^{[\f p3]}729^{p-1}+7p^2U_{p-3}\mod {p^3},
\\&G_{p-1}^{(6)}\e (-1)^{\f{p-1}2}186624^{p-1}+\f{155}9p^2E_{p-3}\mod
{p^3},
\\&G_{p-1}^{(4)}\e (-1)^{[\f p4]}4096^{p-1}+13p^2s_{p-3}\mod {p^3},\endalign$$
where $\{s_n\}$ is given by $s_0=1$ and $s_n=1-\sum_{k=0}^{n-1} \b
nk2^{2n-1-2k}s_k\q(n\ge 1)$.
\endpro
Proof. By [S2, Lemma 2.9], for $k=1,2,\ldots,p-1$,
$$\b{p-1}k(-1)^k\e 1-pH_k+\f{p^2}2\big(H_k^2-H_k^{(2)}\big)\mod
{p^3}.$$ Since $\b{-1-x}k=(-1)^k\b{x+k}k$, appealing to (2.8),
$$\align G_{p-1}(x)&=\sum_{k=0}^{p-1}\b{p-1}k(-1)^k\b xk\b{-1-x}k
\\&\e\sum_{k=0}^{p-1}\b xk\b{-1-x}k-p \sum_{k=1}^{p-1}\b xk\b
{-1-x}kH_k\\&\qq+\f{p^2}2\sum_{k=0}^{p-1}\b
xk\b{x+k}k(-1)^k\big(H_k^2-H_k^{(2)}\big)
 \mod {p^3}.\endalign$$ By (2.13),
$$\sum_{k=0}^{p-1}\b xk\b{-1-x}k\e (-1)^{\xp}+p^2x'(x'+1)E_{p-3}(-x)
\mod {p^3}.$$
 It is well known that $\sum_{r=1}^n\b nr(-1)^r\f 1r=-H_n$. Thus,
 appealing to  (2.7) and (2.17),
 $$\align -(-1)^{\xp}\sum_{k=1}^{p-1}\b xk\b {-1-x}kH_k&\e \sum_{k=1}^{p-1}\b
 xk\b{-1-x}k\f 1k\e -2H_{\xp}+2px'H_{\xp}^{(2)}
 \\&\e -2\f{B_{p^2(p-1)}(-x)-B_{p^2(p-1)}}{p^2(p-1)}
 \mod {p^2}.\endalign$$
 By [T2],
 $$\sum_{k=1}^n\b nk\b{n+k}k(-1)^kH_k^{(2)}=-2(-1)^n\sum_{k=1}^n
 \f{(-1)^k}{k^2}.\tag 3.6$$
 By [LN],
$$\sum_{k=1}^n\b nk\b{n+k}k(-1)^kH_k^2=4(-1)^nH_n^2+2(-1)^n
\sum_{k=1}^n  \f{(-1)^k}{k^2}.\tag 3.7$$ Thus,
$$\align&\sum_{k=0}^{p-1}\b
xk\b{x+k}k(-1)^k\big(H_k^2-H_k^{(2)}\big)
\\&\e \sum_{k=0}^{\xp}\b
{\xp}k\b{\xp+k}k(-1)^k\big(H_k^2-H_k^{(2)}\big) =4(-1)^{\xp}
\Big(H_{\xp}^2+\sum_{k=1}^{\xp}\f{(-1)^k}{k^2} \Big)\mod
p.\endalign$$ By (2.12),
$$\sum_{k=1}^{\xp}\f{(-1)^k}{k^2}\e\f 12(-1)^{\xp}E_{p-3}(-x)\mod
p.$$
 Hence
$$\align G_{p-1}(x)&\e (-1)^{\xp}+p^2x'(x'+1)E_{p-3}(-x)
-p(-1)^{\xp}2\f{B_{p^2(p-1)}(-x)-B_{p^2(p-1)}}{p^2(p-1)}\\&\q+\f{p^2}2
\cdot 4(-1)^{\xp}
\Big(\Ls{B_{p^2(p-1)}(-x)-B_{p^2(p-1)}}{p^2(p-1)}^2+\f
12(-1)^{\xp}E_{p-3}(-x)\Big)\mod{p^3}.\endalign$$ This yields the
result for $G_{p-1}(x)\mod {p^3}$. Now assume that $p>3$. For
$m=3,4,6$ we see that $(-1)^{\langle -\f 1m\rangle_p}=(-1)^{[\f
pm]}$. Taking $x=-\f 13,-\f 14,-\f 16$ in the congruence for
$G_{p-1}(x)\mod {p^3}$ and applying (2.20)-(2.22) and (2.14)-(2.15)
we deduce that
$$\align& G_{p-1}^{(3)}=27^{p-1}G_{p-1}\Big(-\f 13\Big) \e
27^{p-1}(-1)^{[\f p3]}\Big(1-2p\Big(-\f 32\qp 3+\f 34p\qp
3^2\Big)\\&\qq\qq\qq\qq\qq\qq+2p^2\Big(-\f 32\qp
3\Big)^2\Big)+p^2\Big(-\f 13\cdot \f 23+1\Big)\cdot 9U_{p-3}
\\&\qq\e(-1)^{[\f p3]}27^{p-1}(1+p\qp 3)^3+7p^2U_{p-3}
=(-1)^{[\f p3]}729^{p-1}+7p^2U_{p-3}\mod {p^3},
\\&G_{p-1}^{(4)}=64^{p-1}G_{p-1}\Big(-\f 14\Big) \e
64^{p-1}(-1)^{[\f p4]}\Big(1-2p\Big(-3\qp 2+\f 32p\qp
2^2\Big)\\&\qq\qq\qq\qq\qq\qq+2p^2(-3\qp 2)^2\Big)+p^2\Big(-\f
14\cdot \f 34+1\Big)\cdot 16s_{p-3}
\\&\qq\e(-1)^{[\f p4]}64^{p-1}(1+p\qp 2)^6+13p^2s_{p-3}
=(-1)^{[\f p4]}4096^{p-1}+13p^2s_{p-3}\mod {p^3},
\\&G_{p-1}^{(6)}=432^{p-1}G_{p-1}\Big(-\f 16\Big)
\\&\qq \e 432^{p-1}(-1)^{[\f p6]}\Big(1-2p\Big(-2\qp 2-\f 32\qp
3+p\Big(\qp 2^2+\f 34\qp
3^2\Big)\Big)\\&\qq\qq\qq\qq\qq\qq+2p^2\Big(2\qp 2+\f 32\qp
3\Big)^2\Big)+p^2\Big(-\f 16\cdot \f 56+1\Big)\cdot 20E_{p-3}
\\&\qq\e(-1)^{\f{p-1}2}16^{p-1}(1+p\qp 2)^4\cdot 27^{p-1}(1+p\qp 3)^3
+\f{155}9p^2E_{p-3}
\\&\qq=(-1)^{\f{p-1}2}186624^{p-1}+\f{155}9p^2E_{p-3}\mod {p^3}.
\endalign$$
This proves the theorem.
\par\q
\par{\bf Remark 3.2} In [S20], the author conjectured that
for any prime $p>3$,
$$G_{p-1}\e (-1)^{\f{p-1}2}256^{p-1}+3p^2E_{p-3}\mod{p^3}.$$
This was solved by Liu and Ni [LN], and can be easily deduced from
Theorem 3.4 (with $x=-\f 12$).

\pro{Theorem 3.5} Let $p$ be an odd prime, $x\in\Bbb Z_p$ and
$x\not\e 0\mod p$. Then
$$\align G_p(x)&\e 1-\b{x}p\b{-1-x}p+
2p\f{B_{p^2(p-1)}(-x)-B_{p^2(p-1)}}{p^2(p-1)}
\\&\qq+2p^2\Big(\f{B_{p^2(p-1)}(-x)-B_{p^2(p-1)}}{p^2(p-1)}\Big)^2
+p^2(-1)^{\xp}E_{p-3}(-x)\mod {p^3}.
\endalign$$
Hence, for $p>3$,
$$\align &G_p\e 12+64(-1)^{\f{p-1}2}p^2E_{p-3}\mod {p^3},
\\&G_p^{(3)}\e 21+243(-1)^{[\f p3]}p^2U_{p-3}\mod {p^3},
\\&G_p^{(4)}\e 52+1024(-1)^{[\f p4]}p^2s_{p-3}\mod {p^3},
\\&G_p^{(6)}\e 372+8640(-1)^{\f{p-1}2}p^2E_{p-3}\mod {p^3},
\endalign$$
where $\{s_n\}$ is given by (2.16).
\endpro
Proof. Since $\b{p-1}k(-1)^k\e 1-pH_k\mod {p^2}$,
$$\align
G_p(x)&=\sum_{k=0}^p\b pk(-1)^k\b xk\b{-1-x}k \\&=1-\b xp\b{-1-x}p
-\sum_{k=1}^{p-1}\f pk\b{p-1}{k-1}(-1)^{k-1}\b xk\b {-1-x}k
\\&\e 1-\b xp\b{-1-x}p
-p\sum_{k=1}^{p-1}\f{\b xk\b{-1-x}k}k(1-pH_{k-1})
\\&=1-\b xp\b{-1-x}p
-p\sum_{k=1}^{p-1}\f{\b xk\b{-1-x}k}k\\&\qq+p^2\sum_{k=1}^{p-1} \b
xk\b{-1-x}k\Big(\f{H_k}k-\f 1{k^2}\Big)\mod {p^3}.
\endalign$$
By [T1],
$$\sum_{k=1}^{p-1}\f{\b xk\b{-1-x}k}{k^2}\e -\f 12
\Big(\sum_{k=1}^{p-1}\f{\b xk\b{-1-x}k}k\Big)^2\mod p.$$ By [S10,
Theorem 2.1],
$$\sum_{k=1}^{p-1}\f{\b xk\b{-1-x}k}k\e
 -2\f{B_{p^2(p-1)}(-x)-B_{p^2(p-1)}}{p^2(p-1)}
\mod {p^2}.$$ Using the summation package Sigma in Mathematica, one
can find and prove the identity
$$\sum_{k=1}^n\b nk\b{n+k}k\f{(-1)^k}kH_k
=2\sum_{k=1}^n\f{(-1)^k}{k^2}.\tag 3.8$$ Thus, appealing to (2.12),
$$\align &\sum_{k=1}^{p-1}\b xk\b{-1-x}k\f{H_k}k
\\&=\sum_{k=1}^{p-1}\b xk\b{x+k}k(-1)^k\f{H_k}k \e \sum_{k=1}^{\xp} \b
{\xp}k\b{\xp+k}k(-1)^k\f{H_k}k
\\&=2\sum_{k=1}^{\xp}\f{(-1)^k}{k^2}\e (-1)^{\xp}E_{p-3}(-x)\mod p.
\endalign$$
Now, combining the above gives
$$\align G_p(x)&\e 1-\b{x}p\b{-1-x}p+
2p\f{B_{p^2(p-1)}(-x)-B_{p^2(p-1)}}{p^2(p-1)}
\\&\qq+2p^2\Big(\f{B_{p^2(p-1)}(-x)-B_{p^2(p-1)}}{p^2(p-1)}\Big)^2
+p^2(-1)^{\xp}E_{p-3}(-x)\mod {p^3}.
\endalign$$
Hence, for $p>3$, taking $x=-\f 12,-\f 13,-\f 14,-\f 16$ in the
above congruence for $G_p(x)\mod {p^3}$ and then applying
(2.20)-(2.22), (2.14-2.15) and (2.9) yields
$$\align G_p&=16^pG_p\Big(-\f 12\Big)\e
16^p\Big(1-\b{-1/2}p^2+2p \big(-2\qp 2+p\qp 2^2\big)
\\&\qq+2p^2(-2\qp 2)^2 +p^2(-1)^{\f{p-1}2}E_{p-3}\Ls 12\Big)
\\&=16^p\Big(1-\b{2p}p^216^{-p}-4p\qp 2+10p^2\qp
2^2+p^2(-1)^{\f{p-1}2}\cdot 4E_{p-3}\Big)
\\&\e -4+16(1+p\qp 2)^4(1-4p\qp 2+10p^2\qp 2^2)+16^p\cdot 4
(-1)^{\f{p-1}2}p^2E_{p-3}
\\&\e 12+64(-1)^{\f{p-1}2}p^2E_{p-3}\mod {p^3},
\endalign$$
$$\align G_p^{(3)}&=27^pG_p\Big(-\f 13\Big)
\e 27^p\Big(1-\b{-1/3}p\b{-2/3}p+2p\Big(-\f 32\qp 3+\f 34p\qp
3^2\Big)
\\&\qq\qq\qq\qq+2p^2\Big(-\f 32\qp 3\Big)^2+p^2(-1)^{[\f p3]}E_{p-3}
\Ls 13\Big)
\\&\e 27^p\Big(1-\b{2p}p\b{3p}p27^{-p}-3p\qp 3+6p^2\qp
3^2+9(-1)^{[\f p3]}p^2U_{p-3}\Big)
\\&\e -\b 21\b 31+27(1+p\qp 3)^3(1-3p\qp 3+6p^2\qp 3^2)+27^p
\cdot 9(-1)^{[\f p3]}p^2U_{p-3}
\\&\e -6+27(1+3p\qp 3+3p^2\qp 3^2)(1-3p\qp 3+6p^2\qp 3^2)+243
(-1)^{[\f p3]}p^2U_{p-3}
\\&\e 21+243(-1)^{[\f p3]}p^2U_{p-3}\mod {p^3},
\endalign$$
$$\align G_p^{(4)}&=64^pG_p\Big(-\f 14\Big)
\e 64^p\Big(1-\b{-1/4}p\b{-3/4}p+2p\Big(-3\qp 2+\f 32p\qp 2^2\Big)
\\&\qq\qq\qq\qq+2p^2(-3\qp 2)^2+p^2(-1)^{[\f p4]}E_{p-3}
\Ls 14\Big)
\\&\e -\b{2p}p\b{4p}{2p}+64(1+p\qp 2)^6(1-6p\qp 2+21p^2\qp 2^2)
+64^p\cdot 16(-1)^{[\f p4]}p^2s_{p-3}
\\&\e -\b 21\b 42+64(1+6p\qp 2+15p^2\qp 2^2)
(1-6p\qp 2+21p^2\qp 2^2) \\&\qq+1024(-1)^{[\f p4]}p^2s_{p-3}
\\&\e 52+1024(-1)^{[\f p4]}p^2s_{p-3}\mod {p^3}
\endalign$$ and
$$\align G_p^{(6)}&=432^pG_p\Big(-\f 16\Big)
\e 432^p\Big(1-\b{-1/6}p\b{-5/6}p+2p\Big(-2\qp 2-\f 32\qp
3\\&\qq+p\Big(\qp 2^2+\f 34\qp 3^2\Big)\Big) +2p^2\Big(-2\qp 2-\f
32\qp 3\Big)^2+p^2(-1)^{[\f p6]}E_{p-3} \Ls 16\Big)
\\&\e 432^p\Big(1-\b{3p}p\b{6p}{3p}432^{-p}
-p(4\qp 2+3\qp 3)\\&\qq+p^2(10\qp 2^2+6\qp 3^2+12\qp 2\qp 3)\Big)
+432^p\cdot 20 (-1)^{\f {p-1}2}p^2E_{p-3}
\\&\e -\b 31\b 63+432(1+p\qp 2)^4(1+p\qp 3)^3
\Big(1-p(4\qp 2+3\qp 3)\\&\qq+p^2(10\qp 2^2+6\qp 3^2+12\qp 2\qp
3)\Big)  +432\cdot 20(-1)^{\f {p-1}2}p^2E_{p-3}
\\&\e -60+432(1+4p\qp 2+6p^2\qp 2^2)(1+3p\qp 3+3p^2\qp 3^2)
\\&\qq\times\Big(1-p(4\qp 2+3\qp 3)+p^2(10\qp 2^2+6\qp 3^2+12\qp 2\qp 3)\Big)
\\&\qq +8640(-1)^{\f {p-1}2}p^2E_{p-3}
\\&\e 372+8640(-1)^{\f {p-1}2}p^2E_{p-3}\mod {p^3}.
\endalign$$
This completes the proof.

\pro{Theorem 3.6} Let $p$ be an odd prime, $x\in\Bbb Z_p$ and
$x\not\e 0\mod p$. Then
$$G_{\f{p-1}2}(x)\e \cases \b{\xp}{\xp /2}^2\f 1{4^{\xp}}\mod p
&\t{if $\xp$ is even,}
\\0\mod p&\t{if $\xp$ is odd.}
\endcases$$
Moreover,
$$G_{\f{p-1}2}(x)\e
\cases \sum_{k=0}^{(p-1)/2}\b xk\b{-1-x}k\f{\b{2k}k}{4^k}\mod {p^2}
\q\t{if $\xp$ is even,}
\\-\sum_{k=(p+1)/2}^{\xp}\b{\xp}k\b{\xp+k}k\f{\b{2k}k}{(-4)^k}
\\\q-p\sum_{k=0}^{\xp}\b {\xp}k\b{\xp+k}k\f{\b{2k}k}{(-4)^k}H_k\mod {p^2}
\q\t{if $\xp$ is odd.}
\endcases$$
Hence, for $p>3$,
$$\align &G_{\f{p-1}2}^{(3)}\e \cases 27^{\f{p-1}2}(4x^2-2p)\mod {p^2}
&\t{if $p\e
1\mod 3$  and so $p=x^2+3y^2$,}\\0\mod {p}&\t{if $p\e 2\mod 3$,}
\endcases
\\&G_{\f{p-1}2}^{(4)}\e \cases 8^{p-1}\cdot 4x^2-2p\mod {p^2}&\t{if $p\e
1,3\mod 8$  and so $p=x^2+2y^2$,}\\0\mod {p}&\t{if $p\e 5,7\mod 8$,}
\endcases
\\&G_{\f{p-1}2}^{(6)}\e \cases 432^{\f{p-1}2}\ls p3\cdot 4x^2-2p\mod {p^2}&\t{if $p\e
1\mod 4$  and so $p=x^2+4y^2$,}\\0\mod {p}&\t{if $p\e 3\mod 4$.}
\endcases
\endalign$$
\endpro
Proof. By [S3, Lemma 2.4], for $k=1,2,\ldots,\f{p-1}2$,
$$\b{\f{p-1}2}k(-1)^k\e \f{\b{2k}k}{4^k}\Big(1-p\sum_{i=1}^k
\f 1{2i-1}\Big)=\f{\b{2k}k}{4^k}\Big(1-p\Big(H_{2k}-\f
12H_k\Big)\Big) \mod {p^2}.$$ Thus,
$$ \align G_{\f{p-1}2}(x)&=\sum_{k=0}^{(p-1)/2}
\b{\f{p-1}2}k(-1)^k\b xk\b{-1-x}k\\&\e \sum_{k=0}^{(p-1)/2}
\f{\b{2k}k}{4^k}\b xk\b{-1-x}k -p\sum_{k=0}^{(p-1)/2}
\f{\b{2k}k}{4^k}\b xk\b{-1-x}k\Big(H_{2k}-\f 12H_k\Big)
\\&\e\sum_{k=0}^{(p-1)/2}
\f{\b{2k}k}{(-4)^k}\b xk\b{x+k}k -p\sum_{k=0}^{p-1}
\f{\b{2k}k}{(-4)^k}\b xk\b{x+k}k\Big(H_{2k}-\f 12H_k\Big)
 \mod
{p^2}.\endalign$$ Using summation package Sigma one can find and
prove
$$\sum_{k=0}^n\f{\b{2k}k}{(-4)^k}\b nk\b{n+k}k
\Big(H_{2k}-\f {2-(-1)^n}2H_k\Big) =0.\tag 3.9$$ When $n$ is even,
this identity has been given by Tauraso [T2].
 Thus,
$$\align &\sum_{k=0}^{p-1}
\f{\b{2k}k}{(-4)^k}\b xk\b{x+k}k\Big(H_{2k}-\f 12H_k\Big)
\\&\e
\sum_{k=0}^{\xp} \f{\b{2k}k}{(-4)^k}\b
{\xp}k\b{\xp+k}k\Big(H_{2k}-\f 12H_k\Big)\\& =\cases 0\mod p&\t{if
$\xp$ is even,}
\\\sum_{k=0}^{\xp}\f{\b{2k}k}{(-4)^k}\b {\xp}k\b{\xp+k}kH_k\mod p
&\t{if $\xp$ is odd.}\endcases\endalign$$ When $\xp$ is odd, from
[S8, Theorem 2.5], $$\sum_{k=0}^{p-1} \f{\b{2k}k}{4^k}\b
xk\b{-1-x}k\e 0\mod {p^2}.$$ For $\f p2<k<p$ we have $p\mid
\b{2k}k$. Thus, for $\xp<\f p2$,
$$\sum_{k=0}^{(p-1)/2} \f{\b{2k}k}{4^k}\b
xk\b{-1-x}k\e \sum_{k=0}^{p-1} \f{\b{2k}k}{4^k}\b xk\b{-1-x}k\e
0\mod {p^2}.$$ If $\xp$ is odd and $\xp>\f p2$, then
$$\align &\sum_{k=0}^{(p-1)/2} \f{\b{2k}k}{4^k}\b
xk\b{-1-x}k \\&=\sum_{k=0}^{p-1} \f{\b{2k}k}{4^k}\b xk\b{-1-x}k
-\sum_{k=(p+1)/2}^{p-1}\f{\b{2k}k}{4^k}\b xk\b{-1-x}k
\\&\e-\sum_{k=(p+1)/2}^{p-1}\f{\b{2k}k}{(-4)^k}\b xk\b{x+k}k
\e -\sum_{k=(p+1)/2}^{p-1}\f{\b{2k}k}{(-4)^k}\b {\xp}k\b{\xp+k}k
\mod{p^2}.\endalign$$ Combining the above gives the congruence for
$G_{\f{p-1}2}(x)\mod {p^2}$. By [G,(6.35)], we have the following
identity due to Bell:
$$\sum_{k=0}^n\b nk\b{n+k}k\f{\b{2k}k}{(-4)^k}
=\cases \b{n}{n/2}^2\f 1{4^n}&\t{if $n$ is even,}
\\0&\t{if $n$ is odd.}\endcases\tag 3.10$$
Hence, from the above,
$$\align G_{\f{p-1}2}(x)&\e
\sum_{k=0}^{(p-1)/2} \f{\b{2k}k}{(-4)^k}\b xk\b{x+k}k \e
\sum_{k=0}^{p-1} \f{\b{2k}k}{(-4)^k}\b xk\b{x+k}k \\&\e
\sum_{k=0}^{\xp} \f{\b{2k}k}{(-4)^k}\b {\xp}k\b{\xp+k}k
\\&=\cases \b{\xp}{\xp /2}^2\f 1{4^{\xp}}\mod p
&\t{if $\xp$ is even,}
\\0\mod p&\t{if $\xp$ is odd.}
\endcases\endalign$$
From [M] or [Su3] we know that
$$\align &\sum_{k=0}^{p-1}\f{\b{2k}k^2\b{3k}k}{108^k}\e\cases
4x^2-2p\mod {p^2}&\t{if $p\e 1\mod 3$  and so $p=x^2+3y^2$,}\\0\mod
{p}&\t{if $p\e 2\mod 3$,}
\endcases
\\&\sum_{k=0}^{p-1}\f{\b{2k}k^2\b{4k}{2k}}{256^k}\e
 \cases 4x^2-2p\mod {p^2}&\t{if $p\e
1,3\mod 8$  and so $p=x^2+2y^2$,}\\0\mod {p}&\t{if $p\e 5,7\mod 8$,}
\endcases
\\&\sum_{k=0}^{p-1}\f{\b{2k}k\b{3k}k\b{6k}{3k}}{1728^k}\e
\cases \ls p3(4x^2-2p)\mod {p^2}&\t{if $p\e 1\mod 4$ and so
$p=x^2+4y^2$,}\\0\mod {p}&\t{if $p\e 3\mod 4$.}
\endcases
\endalign$$
For $p>3$, taking $x=-\f 13,-\f 14,-\f 16$ in the congruence for
$G_{\f{p-1}2}(x)\mod {p^2}$ yields the congruences for $
G_{\f{p-1}2}^{(3)},\ G_{\f{p-1}2}^{(4)}$ and $G_{\f{p-1}2}^{(6)}$.
\par\q\par{\bf Remark 3.3} In [S20], the author proved that for any
prime $p>3$,
$$G_{\f{p-1}2}\e \cases 4^px^2-2p\mod {p^2}&
\t{if $p=x^2+4y^2\e 1\mod 4$,}
\\0\mod {p^2}&\t{if $p\e 3\mod 4$.}
\endcases$$

 \pro{Theorem 3.7} Let $p$ be an odd prime and $m,x\in\Bbb
Z_p$. Then
$$\align &\sum_{k=0}^{p-1}\f{G_k(x)}{m^k}\e
\sum_{k=0}^{p-1}\f{\b xk\b{-1-x}k}{(1-m)^k}\mod p\qtq{for}m\not\e
0,1\mod p,
\\&\sum_{k=0}^{p-1}\b{2k}k\f{G_k(x)}{(m+2)^k}\\&\qq\e
\Ls{(m+2)(m-2)}p\sum_{k=0}^{p-1}\b{2k}k\f{\b xk\b{-1-x}k}{(2-m)^k}
\mod p\q \t{for}\q m\not\e \pm 2\mod p.\endalign$$ Hence, for $p>3$,
$$\align &\sum_{k=0}^{p-1}\f{G_k^{(3)}}{m^k}
\e \sum_{k=0}^{p-1}\f{\b{2k}k\b{3k}k}{(27-m)^k}\mod
p\qtq{for}m\not\e 0,27\mod p,
\\ &\sum_{k=0}^{p-1}\f{G_k^{(4)}}{m^k}
\e \sum_{k=0}^{p-1}\f{\b{2k}k\b{4k}{2k}}{(64-m)^k}\mod p
\qtq{for}m\not\e 0,64\mod p,
\\ &\sum_{k=0}^{p-1}\f{G_k^{(6)}}{m^k}
\e \sum_{k=0}^{p-1}\f{\b{3k}k\b{6k}{3k}}{(432-m)^k}\mod p
\qtq{for}m\not\e 0,432\mod p,
\\&\sum_{k=0}^{p-1}\b{2k}k\f{G_k^{(3)}}{m^k}
\e\Ls {m(m-108)}p\sum_{k=0}^{p-1}\f{\b{2k}k^2\b{3k}k}{(108-m)^k}\mod
p\qtq{for}m\not\e 0,108\mod p,
\\&\sum_{k=0}^{p-1}\b{2k}k\f{G_k^{(4)}}{m^k}
\e\Ls
{m(m-256)}p\sum_{k=0}^{p-1}\f{\b{2k}k^2\b{4k}{2k}}{(256-m)^k}\mod p
\qtq{for}m\not\e 0,256\mod p,
\\&\sum_{k=0}^{p-1}\b{2k}k\f{G_k^{(6)}}{m^k}
\e\Ls {m(m-1728)}p\sum_{k=0}^{p-1}\f{\b{2k}k\b{3k}k\b{6k}{3k}}
{(1728-m)^k}\mod p\ \t{for}\ m\not\e 0,1728\mod p.
\endalign$$
\endpro
Proof. Taking $u_k=(-1)^k\b xk\b{-1-x}k$ and $v_k=G_k(x)$ in Lemma
2.1 yields the congruence for  $\sum_{k=0}^{p-1}\f{G_k(x)}{m^k}$
$\mod p$. Taking $x=-\f 13,-\f 14,-\f 16$ and then replacing $m$
with $\f m{27},\f m{64},\f m{432}$ yields the congruences for
$\sum_{k=0}^{p-1}\f{G_k^{(r)}}{m^k}\mod p$ $(r=3,4,6)$.
\par Note that $p\mid \b{2k}k$ for $\f p2<k<p$. Taking
 $u_k=(-1)^k\b xk\b{-1-x}k$ and $v_k=G_k(x)$ in Lemma 2.2 yields the
congruence for  $\sum_{k=0}^{p-1}\b{2k}k\f{G_k(x)}{(m+2)^k}\mod p$.
Taking $x=-\f 13,-\f 14,-\f 16$ and then replacing $m$ with $\f
m{27}-2,\f m{64}-2,\f m{432}-2$ yields the congruences for
$\sum_{k=0}^{p-1} \b{2k}k\f{G_k^{(r)}}{m^k}\mod p$ $(r=3,4,6)$.
Thus, the theorem is proved.

\pro{Theorem 3.8} Let $p>3$ be a prime. Then
$$\align &\sum_{n=0}^{p-1}\f{G_n^{(3)}}{(-27)^n}
\e\cases 2x\mod{p}&\t{if $p=x^2+3y^2\e 1\mod 3$ and $3\mid x-1$,}
\\0\mod p
&\t{if $p\e 2\mod 3$,}\endcases
\\&\sum_{n=0}^{p-1}\f{G_n^{(3)}}{3^n}\e
\sum_{n=0}^{p-1}\f{G_n^{(3)}}{243^n}
\e\cases -L\mod{p}&\t{if $3\mid p-1$, $4p=L^2+27M^2$ and $3\mid
L-1$,}
\\0\mod p&\t{if $p\e 2\mod 3$.}\endcases
\endalign$$
\endpro
Proof. By Theorem 3.7, $$\sum_{k=0}^{p-1}\f{G_k^{(3)}}{m^k} \e
\sum_{k=0}^{p-1}\f{\b{2k}k\b{3k}k}{(27-m)^k}\mod p.$$ Now, taking
$m=-27,3,243$ and then applying [S8, Theorem 3.4] and [S6, Theorem
3.2] yields the result.

\pro{Theorem 3.9} Let $p>3$ be a prime. Then
$$\align &\Ls{-3}p\sum_{n=0}^{p-1}G_n^{(4)}\e
\Ls 6p\sum_{n=0}^{p-1}\f{G_n^{(4)}}{4096^n}
 \\&\qq\e\cases 2x\mod{p}&\t{if
$p=x^2+7y^2\e 1,2,4\mod 7$ and $\sls x7=1$,}
\\0\mod p
&\t{if $p\e 3,5,6\mod 7$,}\endcases
\\&\Ls 6p\sum_{n=0}^{p-1}\f{G_n^{(4)}}{(-8)^n}
\e \Ls 3p\sum_{n=0}^{p-1}\f{G_n^{(4)}}{(-512)^n}
\\&\qq\e\cases
2x\mod{p}&\t{if $p=x^2+4y^2\e 1\mod 4$ and $4\mid x-1$,}
\\0\mod p&\t{if $p\e 3\mod 4$,}\endcases
\\&\sum_{n=0}^{p-1}\f{G_n^{(4)}}{16^n}
\e\Ls{-2}p\sum_{n=0}^{p-1}\f{G_n^{(4)}}{256^n}
 \\&\qq\e\cases 2x\mod{p}&\t{if $p=x^2+3y^2\e 1\mod 3$
and $3\mid x-1$,}
\\0\mod {p}
&\t{if $p\e 2\mod 3$,}
\endcases
\\&\sum_{n=0}^{p-1}\f{G_n^{(4)}}{(-64)^n}
 \e\cases (-1)^{[\f p8]+\f{p-1}2}2x
 \mod{p}&\t{if $p=x^2+2y^2\e 1,3\mod 8$
and $4\mid x-1$,}
\\0\mod {p}
&\t{if $p\e 5,7\mod 8$.}
\endcases
\endalign$$
\endpro
Proof. From Theorem 3.7,
$$\sum_{k=0}^{p-1}\f{G_k^{(4)}}{m^k} \e
\sum_{k=0}^{p-1}\f{\b{2k}k\b{4k}{2k}}{(64-m)^k}\mod p.$$ Now taking
$m=1,4096,-8,-512,16,256,-64$ and then applying [S5, Theorems
2.2-2.4 and 4.3] yields the result.
\par\q
\par Based on calculations with Maple,  we pose the following
conjectures.
 \pro{Conjecture 3.1} Let $p>3$ be a prime. Then
$$\align &\sum_{n=0}^{p-1}\f{G_n^{(3)}}{(-27)^n}
\e\cases 2x-\f p{2x}\mod{p^2}&\t{if $p=x^2+3y^2\e 1\mod 3$ and
$3\mid x-1$,}
\\-\f 32p\b{\f{p-1}2}{\f{p-5}6}^{-1}\mod {p^2}
&\t{if $p\e 2\mod 3$,}\endcases
\\&\sum_{n=0}^{p-1}\f{G_n^{(3)}}{3^n}
\e\cases -L+\f pL\mod{p^2}&\t{if $3\mid p-1$, $4p=L^2+27M^2$ and
$3\mid L-1$,}
\\-\f 43p\big(\f{p-2}3!\big)^3\mod {p^2}
&\t{if $p\e 2\mod 3$,}\endcases
\\&\sum_{n=0}^{p-1}\f{G_n^{(3)}}{243^n}
\e\cases -L+\f pL\mod{p^2}&\t{if $3\mid p-1$, $4p=L^2+27M^2$ and
$3\mid L-1$,}
\\0\mod {p^2}
&\t{if $p\e 2\mod 3$.}
\endcases\endalign$$
\endpro

\par\q
\pro{Conjecture 3.2} Let $p>3$ be a prime. Then
$$\align&\Ls p3\sum_{n=0}^{p-1}G_n^{(4)}
\e \f{21\sls p7-19}2\Ls 6p\sum_{n=0}^{p-1}\f{G_n^{(4)}}{4096^n}
\\&\qq\e\cases 2x-\f p{2x}\mod{p^2}&\t{if $p=x^2+7y^2\e 1,2,4\mod 7$ with
$\ls x7=1$,}
\\-\f{5p}{\b{3[p/7]}{[p/7]}}\mod{p^2}&\t{if $p\e 3\mod 7$,}
\\\f{15p}{4\b{3[p/7]}{[p/7]}}\mod{p^2}&\t{if $p\e 5\mod 7$,}
\\-\f{25p}{22\b{3[p/7]}{[p/7]}}\mod{p^2}&\t{if $p\e 6\mod 7$,}
\endcases
\\&\Ls 6p\sum_{n=0}^{p-1}\f{G_n^{(4)}}{(-8)^n}
\e \big(3-2(-1)^{\f{p-1}2}\big)\Ls
3p\sum_{n=0}^{p-1}\f{G_n^{(4)}}{(-512)^n}
\\&\qq\e\cases 2x-\f
p{2x}\mod{p^2}&\t{if $p=x^2+4y^2\e 1\mod 4$ with $4\mid x-1$,}
\\ \f{5p}2\b{\f{p-3}2}{\f{p-3}4}^{-1}\mod {p^2}
&\t{if $p\e 3\mod 4$,}\endcases
\\&\sum_{n=0}^{p-1}\f{G_n^{(4)}}{16^n}\e
\Big(4-3\Ls{-3}p\Big)\Ls{-2}p \sum_{n=0}^{p-1}\f{G_n^{(4)}}{256^n}
\\&\qq
 \e\cases 2x-\f p{2x}\mod{p^2}&\t{if $p=x^2+3y^2\e 1\mod 3$ with
$3\mid x-1$,}
\\ -\f 72p\b{\f{p-1}2}{\f{p-5}6}^{-1}\mod {p^2}
&\t{if $p\e 2\mod 3$,}\endcases
\\&\sum_{n=0}^{p-1}\f{G_n^{(4)}}{(-64)^n} \e\cases
(-1)^{[\f p8]+\f{p-1}2}(2x-\f p{2x})\mod{p^2} \\\qq\qq\qq\t{if
$p=x^2+2y^2\e 1,3\mod 8$ with $4\mid x-1$,}
\\ -\f {4}{2-(-1)^{\f{p-1}2}}p\b{\f{p-1}2}{[\f p8]}^{-1}\mod {p^2}
\q\t{if $p\e 5,7\mod 8$.}\endcases
\endalign$$
\endpro

\pro{Conjecture 3.3} Let $p$ be a prime with $p\e 1\mod 3$. Then
$$\sum_{n=1}^{p-1}\f{nG_n^{(3)}}{243^n}\e 0\mod {p^2}.$$
\endpro

\pro{Conjecture 3.4} Let $p>3$ be a prime,
$m\in\{-3267,-1350,-108,44,100,135,$
$300,1836,8748,110700,27000108\}$ and $p\nmid m(108-m)$. Then
$$\sum_{k=0}^{p-1}\b{2k}k\f{G_k^{(3)}}{m^k}
\e \Ls{m(m-108)}p\sum_{k=0}^{p-1}\f{\b{2k}k^2\b{3k}k}{(108-m)^k}
\mod {p^2}.$$
\endpro

\pro{Conjecture 3.5} Let $p>3$ be a prime,
$m\in\{-24591257600,-2508800,-614400,$ $-20480,-2048,-392,175,
400,1280,4225,12544,83200,6635776,199148800\}$ and $p\nmid
m(256-m)$. Then
$$\sum_{k=0}^{p-1}\b{2k}k\f{G_k^{(4)}}{m^k}
\e \Ls{m(m-256)}p\sum_{k=0}^{p-1}\f{\b{2k}k^2\b{4k}{2k}}{(256-m)^k}
\mod {p^2}.$$
\endpro

\pro{Conjecture 3.6} Let $p>3$ be a prime, $m\in\{-16579647,-285768,
-52272, -6272,$ $5103,34496,886464,12289728, 884737728,
147197953728, 262537412640769728\}$ and $p\nmid m$ $(1728-m)$. Then
$$\sum_{k=0}^{p-1}\b{2k}k\f{G_k^{(6)}}{m^k}
\e
\Ls{m(m-1728)}p\sum_{k=0}^{p-1}\f{\b{2k}k\b{3k}k\b{6k}{3k}}{(1728-m)^k}
\mod {p^2}.$$
\endpro
\par{\bf Remark 3.4} Let $p$ be a prime with $p>3$. For the values
of $m$ in Conjectures 3.4-3.6, the congruences for
$$\sum_{k=0}^{p-1}\f{\b{2k}k^2\b{3k}k}{(108-m)^k},\q
\sum_{k=0}^{p-1}\f{\b{2k}k^2\b{4k}{2k}}{(256-m)^k},\q
\sum_{k=0}^{p-1}\f{\b{2k}k\b{3k}k\b{6k}{3k}}{(1728-m)^k} \mod
{p^2}$$ were conjectured by Z.W. Sun in [Su1,Su4] and the author in
[S3], and partially solved by the author in [S5-S7].

\pro{Conjecture 3.7} Let $p$ be an odd prime and $m,r\in\Bbb Z^+$.
Then
$$\align G_{mp^r}^{(3)}\e G_{mp^{r-1}}^{(3)}\mod {p^{2r}},
\ G_{mp^r}^{(4)}\e G_{mp^{r-1}}^{(4)}\mod {p^{2r}}, \
G_{mp^r}^{(6)}\e G_{mp^{r-1}}^{(6)}\mod {p^{2r}} .\endalign$$
\endpro

\pro{Conjecture 3.8} Suppose that $p$ is an odd prime,
$m\in\{1,3,5,\ldots\}$ and $r\in\{2,3,4,\ldots\}$. Then
$$\align &
G_{\f{mp^r-1}2}^{(3)}\e p^2G_{\f{mp^{r-2}-1}2}^{(3)}\mod {p^{2r-1}}
\qtq{for}p\e 2\mod 3,
\\&G_{\f{mp^r-1}2}^{(4)}\e p^2G_{\f{mp^{r-2}-1}2}^{(4)}\mod {p^{2r-1}}
\qtq{for}p\e 5,7\mod 8,
\\&G_{\f{mp^r-1}2}^{(6)}\e p^2G_{\f{mp^{r-2}-1}2}^{(6)}\mod {p^{2r-1}}
\qtq{for}p\e 3\mod 4.
\endalign$$
\endpro
\pro{Conjecture 3.9} Suppose that $p$ is an odd prime,
$m\in\{1,3,5,\ldots\}$ and $r\in\{2,3,4,\ldots\}$. Then
$$\align &G_{\f{mp^r-1}2}^{(3)}\e (-1)^{\f{p-1}2}(4x^2-2p)G_{\f{mp^{r-1}-1}2}^{(3)}-
p^2G_{\f{mp^{r-2}-1}2}^{(3)}\mod {p^r}\\&\qq\qq\qq\qq\qq\qq
\qq\qq\qq\ \t{for}\q p=x^2+3y^2\e 1\mod 3,
\\&G_{\f{mp^r-1}2}^{(4)}\e (4x^2-2p)G_{\f{mp^{r-1}-1}2}^{(4)}-
p^2G_{\f{mp^{r-2}-1}2}^{(4)}\mod {p^r} \ \t{for}\ p=x^2+2y^2\e
1,3\mod 8,
\\&G_{\f{mp^r-1}2}^{(6)}\e (4x^2-2p)G_{\f{mp^{r-1}-1}2}^{(6)}-
p^2G_{\f{mp^{r-2}-1}2}^{(6)}\mod {p^r}\ \t{for}\ p=x^2+4y^2\e 1\mod
4.
\endalign$$
\endpro

\pro{Conjecture 3.10} Let $p$ be a prime with $p>3$. Then
$$\align&G_{2p}\e G_2+3072(-1)^{\f{p-1}2}p^2E_{p-3}\mod {p^3},
\\&G_{3p}\e G_3+94464(-1)^{\f{p-1}2}p^2E_{p-3}\mod {p^3},
\\&G_{2p}^{(3)}\e G_2^{(3)}+20412(-1)^{[\f p3]}p^2U_{p-3}\mod {p^3},
\\&G_{2p-1}\e (-1)^{\f{p-1}2}16^{4(p-1)}G_1+164p^2E_{p-3}\mod {p^3},
\\&G_{2p-1}^{(3)}\e (-1)^{[\f p3]}27^{4(p-1)}G_1^{(3)}+660
p^2U_{p-3}\mod {p^3},
\\&G_{2p-1}^{(6)}\e (-1)^{\f{p-1}2}432^{4(p-1)}G_1^{(6)}+\f {82580}3
p^2E_{p-3}\mod {p^3}.
\endalign$$
\endpro

\pro{Conjecture 3.11} Suppose $n\in\Bbb Z^+$. If $x\in(-1,0)$, then
$G_n(x)^2<G_{n+1}(x)G_{n-1}(x)$. If $x\not\in [-1,0]$, then $
G_n(x)^2>G_{n+1}(x)G_{n-1}(x)$.
\endpro

\section*{4. Congruences for $V_n(x)$}
Recall that
$$V_n(x)=\sum_{k=0}^n\b nk\b{n+k}k(-1)^k \b
xk\b{-1-x}k.$$
 Using sumtools in Maple we find that
$$\aligned &(n+1)^3m^{n+1}V_{n+1}(x)
=(2n+1)(mn(n+1)+m(2x^2+2x+1))m^nV_n(x)
\\&\qq\qq\qq\qq\qq\q-m^2n^3m^{n-1}V_{n-1}(x)\
(n\ge 1).\endaligned\tag 4.1$$ Thus,  $m^nV_n(x)$ is an Ap\'ery-like
sequence of the first kind with $a=m$, $b=m(2x^2+2x+1)$ and $c=m^2$.
\pro{Theorem 4.1} For $n=0,1,2,\ldots$ we have
$$V_n(x)=\sum_{k=0}^n\b xk^2\b{-1-x}{n-k}^2
=\sum_{k=0}^n\b nk\b{n+k}k(-1)^{n-k}G_k(x).$$
\endpro
Proof. Set $V'_n(x)=\sum_{k=0}^n\b xk^2\b{-1-x}{n-k}^2$. Then
$V'_0(x)=1=V_0(x)$ and $V'_1(x)=x^2+(x+1)^2=V_1(x)$. Using sumtools
in Maple we find that
$$(n+1)^3V'_{n+1}(x)=(2n+1)(n(n+1)+2x^2+2x+1)V'_n(x)-n^3V'_{n-1}(x)
\q(n\ge 1).$$ Thus $V'_n(x)=V_n(x)$ for $n=0,1,2,\ldots$ by (4.1).
Now applying (2.6),
$$\sum_{k=0}^n\b nk\b{n+k}k(-1)^{n-k}G_k(x)=\sum_{k=0}^n\b nk\b{n+k}k
(-1)^k\b xk\b{-1-x}k=V_n(x).$$ This proves the theorem.
\par\q
\par For $n=0,1,2,\ldots$, let
$V_n^{(3)},V_n^{(4)}$ and $V_n^{(6)}$ be given by (1.11)-(1.13).
Then $V_n^{(3)},V_n^{(4)}$ and $V_n^{(6)}$ are Ap\'ery-like integral
sequences of the first kind with $(a,b,c)=(27,15,729),$
$(64,40,4096)$ and $(432,312,186624)$, respectively.

\pro{Theorem 4.2} Let $p>3$ be a prime and $x\in\Bbb Z_p$. Then
$$\align V_p(x)&\e 1-2\b
xp\b{-1-x}p +2p \f{B_{p^2(p-1)}(-x)-B_{p^2(p-1)}}{p^2(p-1)}
\\&\qq+2p^2\Big(\f{B_{p^2(p-1)}(-x)-B_{p^2(p-1)}}{p^2(p-1)}\Big)^2\mod
{p^3}.\endalign$$ Taking $x=-\f 12,-\f 13,-\f 14,-\f 16$ yields
$$\align V_p\e 8\mod {p^3},\q V_p^{(3)}\e 15\mod {p^3},
\q V_p^{(4)}\e 40\mod {p^3},\q V_p^{(6)}\e 312\mod {p^3}.
\endalign$$
\endpro
Proof. It is clear that
$$\align &V_p(x)-1+\b{2p}p\b
xp\b{-1-x}p
\\&=\sum_{k=1}^{p-1}\b pk\b {p+k}k(-1)^k\b xk\b{-1-x}k
\\&=p\sum_{k=1}^{p-1}\f{(p+k)(p^2-1^2)\cdots(p^2-(k-1)^2)}{k!^2}
(-1)^k\b xk\b{-1-x}k
\\&\e p\sum_{k=1}^{p-1}\f{(p+k)(-1)^{k-1}
(k-1)!^2(1-p^2H_{k-1}^{(2)})}{k!^2} (-1)^k\b xk\b{-1-x}k
\\&\e-p\sum_{k=1}^{p-1}\f{p+k-k\cdot p^2H_{k-1}^{(2)}}{k^2}\b
xk\b{-1-x}k\mod {p^4}.\endalign$$ Thus, appealing to (2.17), (2.18)
and the known congruence $\b{2p}p=2\b{2p-1}{p-1}\e 2\mod {p^3}$,
$$\align &V_p(x)-1+2\b
xp\b{-1-x}p
\\&\e -p^2\sum_{k=1}^{p-1}\f{\b xk\b{-1-x}k}{k^2}-p\sum_{k=1}^{p-1}
\f{\b xk\b{-1-x}k}k
\\&\e 2p^2\Big(\f{B_{p^2(p-1)}(-x)-B_{p^2(p-1)}}{p^2(p-1)}\Big)^2 +2p
\f{B_{p^2(p-1)}(-x)-B_{p^2(p-1)}}{p^2(p-1)}\mod {p^3}.
\endalign$$
Since $\b{mp}{np}\e \b mn\mod{p^3}$, taking $x=-\f 12$ in the above
congruence and then applying (2.19) we get
$$\align V_p&=16^pV_p\Big(-\f 12\Big)
\e 16^p\Big(1-2\b{2p}p^216^{-p}+2p(-2\qp 2+p\qp 2^2)+2p^2(-2\qp
2)^2\Big)
\\&\e -8+16(1+p\qp 2)^4(1-4p\qp 2+10p^2\qp 2^2)
\\&\e -8+16(1+4p\qp 2+6p^2\qp 2^2)(1-4p\qp 2+10p^2\qp 2^2)
\e 8\mod {p^3}.\endalign$$ Taking $x=-\f 13$ in the congruence for
$V_p(x)\mod {p^3}$ and then applying (2.20) we get
$$\align V_p^{(3)}&=27^pV_p\Big(-\f 13\Big)
\\&\e 27^p\Big(1-2\b{2p}p\b{3p}p27^{-p}+2p\Big(-\f 32\qp 3+\f 34p\qp
3^2\Big)+2p^2\Big(-\f 32\qp 3\Big)^2\Big)
\\&\e -2\cdot 2\cdot 3+27(1+p\qp 3)^3(1-3p\qp 3+6p^2\qp 3^2)
\\&\e -12+27(1+3p\qp 3+3p^2\qp 3^2)(1-3p\qp 3+6p^2\qp 3^2)
\e 15\mod {p^3}.\endalign$$ Taking $x=-\f 14$ in the congruence for
$V_p(x)\mod {p^3}$ and then applying (2.21) we get
$$\align V_p^{(4)}&=64^pV_p\Big(-\f 14\Big)
\\&\e 64^p\Big(1-2\b{2p}p\b{4p}{2p}64^{-p}+2p\Big(-3\qp 2+\f 32p\qp
2^2\Big)+2p^2(-3\qp 2)^2\Big)
\\&\e -2\b 21\b 42+64(1+p\qp 2)^6(1-6p\qp 2+21p^2\qp 2^2)
\\&\e -24+64(1+6p\qp 2+15p^2\qp 2^2)(1-6p\qp 2+21p^2\qp 2^2)
\e 40\mod {p^3}.\endalign$$ Taking $x=-\f 16$ in the congruence for
$V_p(x)\mod {p^3}$ and then applying (2.22) yields
$$\align V_p^{(6)}&=432^pV_p\Big(-\f 16\Big)
\\&\e 432^p\Big(1-2\b{3p}p\b{6p}{3p}432^{-p}+2p\Big(-2\qp 2-\f 32\qp
3+p\Big(\qp 2^2+\f 34\qp 3^2\Big)\Big)\\&\qq\qq\qq\qq+2p^2\Big(-2\qp
2-\f 32\qp 3\Big)^2\Big)
\\&\e -2\b 31\b 63+432(1+p\qp 2)^4(1+p\qp 3)^3
\Big(1-p(4\qp 2+3\qp 3)\\&\qq+p^2(10\qp 2^2+6\qp 3^2+12\qp 2\qp
3)\Big)
\\&\e -120+432=312\mod {p^3}.\endalign$$
This completes the proof.

\pro{Theorem 4.3} Let $p$ be an odd prime and $x\in\Bbb Z_p$. Then
$$ V_{p-1}(x)\e 1-2p\f{B_{p^2(p-1)}(-x)-B_{p^2(p-1)}}{p^2(p-1)}
+2p^2\Big(\f{B_{p^2(p-1)}(-x)-B_{p^2(p-1)}}{p^2(p-1)}\Big)^2\mod
{p^3}.$$ Taking $x=-\f 13,-\f 14,-\f 16$ yields that for $p>3$,
$$  V_{p-1}^{(3)}\e 729^{p-1}\mod {p^3},
\q V_{p-1}^{(4)}\e 4096^{p-1}\mod {p^3}, \q V_{p-1}^{(6)}\e
186624^{p-1}\mod {p^3}. $$
\endpro
Proof. Clearly
$$\align &V_{p-1}(x)\\&=\sum_{k=0}^{p-1}\b{p-1}k\b{p-1+k}k(-1)^k\b
xk\b{-1-x}k
\\&=1+\sum_{k=1}^{p-1}\f{p(p-k)(p^2-1^2)(p^2-2^2)\cdots(p^2-(k-1)^2)}{k!^2}
(-1)^k\b xk\b{-1-x}k
\\&\e
1+p\sum_{k=1}^{p-1}\f{(p-k)(-1)^{k-1}(k-1)!^2(1-p^2H_{k-1}^{(2)})}{k!^2}
(-1)^k\b xk\b{-1-x}k
\\&\e 1-p^2\sum_{k=1}^{p-1}\f{\b xk\b{-1-x}k}{k^2}
+p\sum_{k=1}^{p-1}\f{\b xk\b{-1-x}k}k -p^3\sum_{k=1}^{p-1}\b
xk\b{-1-x}k\f{H_{k-1}^{(2)}}k\mod {p^4}.
\endalign$$
Hence, applying (2.17) and (2.18) yields
$$V_{p-1}(x)\e 1-2p\f{B_{p^2(p-1)}(-x)-B_{p^2(p-1)}}{p^2(p-1)}
+2p^2\Big(\f{B_{p^2(p-1)}(-x)-B_{p^2(p-1)}}{p^2(p-1)}\Big)^2\mod
{p^3}.$$ Now assume that $p>3$. Taking $x=-\f 13$ and then applying
(2.20) yields
$$\align V_{p-1}^{(3)}&=27^{p-1}V_{p-1}\Big(-\f 13\Big)
\e 27^{p-1}\Big(1-2p\Big(-\f 32\qp 3+\f 34p\qp 3^2\Big)+2p^2
\Big(-\f 32\qp 3\Big)^2\Big)
\\&\e 27^{p-1}(1+3p\qp 3+3p^2\qp 3^2)\e 27^{p-1}(1+p\qp 3)^3
=729^{p-1}\mod {p^3}.\endalign$$ Taking $x=-\f 14$ and then applying
(2.21) yields
$$\align V_{p-1}^{(4)}&=64^{p-1}V_{p-1}\Big(-\f 14\Big)
\e 64^{p-1}\Big(1-2p\Big(-3\qp 2+\f 32p\qp 2^2\Big)+2p^2 (-3\qp
2)^2\Big)
\\&\e 64^{p-1}(1+6p\qp 2+15p^2\qp 2^2)\e 64^{p-1}(1+p\qp 2)^6
=4096^{p-1}\mod {p^3}.\endalign$$ Taking $x=-\f 16$ and then
applying (2.22) yields
$$\align V_{p-1}^{(6)}&=432^{p-1}V_{p-1}\Big(-\f 16\Big)
\e 432^{p-1}\Big(1-2p\Big(-2\qp 2-\f 32\qp 3+p\Big(\qp 2^2+\f 34\qp
3^2\Big)\Big)\\&\qq+2p^2\Big (-2\qp 2-\f 32\qp 3\Big)^2\Big)
\\&\e (1+p\qp 2)^4(1+p\qp 3)^3(1+4p\qp 2+6p^2\qp 2^2)(1+3p\qp
3+3p^2\qp 3^2)\\&\e 186624^{p-1} \mod {p^3}.
\endalign$$
This completes the proof.
\par \q
\par{\bf Remark 4.1} In [S20], the author conjectured that for any prime
$p>3$, $$V_{p-1}\e 256^{p-1}-\f 32p^3B_{p-3}\mod {p^4},$$ and proved
the congruence modulo $p^3$.

\pro{Theorem 4.4} Let $p>3$ be a prime, $x\in\Bbb Z_p$, $2x+1\not\e
0\mod p$ and $x'=(x-\xp)/p$. Then
$$\sum_{n=0}^{p-1}V_n(x)\e
\f{1+2x'}{1+2x}p+\f{x'(x'+1)+1}{1+2x}p^3B_{p-2}(-x)
 \mod {p^4}.$$ Hence
$$\align &\sum_{n=0}^{p-1}\f{V_n^{(3)}}{27^n}\e (-1)^{[\f p3]}p+14p^3U_{p-3}
\mod {p^4},
\\&\sum_{n=0}^{p-1}\f{V_n^{(4)}}{64^n}\e (-1)^{\f
{p-1}2}p+13p^3E_{p-3} \mod {p^4},
\\&\sum_{n=0}^{p-1}\f{V_n^{(6)}}{432^n}\e (-1)^{[\f p3]}p+\f{155}4p^3U_{p-3}
\mod {p^4} .\endalign$$
\endpro
Proof. By the definition of $V_n(x)$ and (2.2),
$$\align\sum_{n=0}^{p-1}V_n(x)&=\sum_{n=0}^{p-1}
\sum_{k=0}^n\b nk\b{n+k}k(-1)^k\b xk\b{-1-x}k
\\&=\sum_{k=0}^{p-1}\b{x}k\b{-1-x}k(-1)^k\b{2k}k\sum_{n=k}^{p-1}
\b{n+k}{2k}
\\&=\sum_{k=0}^{p-1}\b{x}k\b{-1-x}k(-1)^k\b{2k}k\b{p+k}{2k+1}
\\&=\sum_{k=0}^{p-1}\b{x}k\b{-1-x}k
\f{p(-1)^k(p^2-1^2)\cdots(p^2-k^2)}{(2k+1)\cdot k!^2}
\\&\e \sum_{k=0}^{p-1}\b{x}k\b{-1-x}k\f{p(1-p^2H_k^{(2)})}{2k+1}
\mod {p^4}.
\endalign$$
By [S12, Theorem 2.1],
 $$\sum_{k=0}^{p-1}\b{x}k\b{-1-x}k\f 1{2k+1}\e
 \f{1+2x'}{1+2x}+p^2\f{x'(x'+1)}{1+2x}B_{p-2}(-x)\mod {p^3}.$$
 Recently Liu and Ni [LN] noted the identity
 $$\sum_{k=0}^n\b
 nk\b{n+k}k\f{(-1)^k}{2k+1}H_k^{(2)}=-2\f{H_n^{(2)}}{2n+1}.\tag 4.2$$
 Since $H_{\f{p-1}2}^{(2)}=\sum_{r=1}^{(p-1)/2}\f 1{r^2}\e 0\mod p$
by [S1, Theorem 5.2], $\f{H_k^{(2)}}{2k+1}\in\Bbb Z_p$ for
$k=0,1,\ldots,p-1$. Now, from the above and (2.11),
$$\align&\sum_{k=0}^{p-1}\b xk\b{-1-x}k\f{H_k^{(2)}}{2k+1}
\\&=\sum_{k=0}^{p-1}\b xk\b{x+k}k(-1)^k\f{H_k^{(2)}}{2k+1}
\e \sum_{k=0}^{\xp}\b {\xp}k\b{\xp+k}k(-1)^k\f{H_k^{(2)}}{2k+1}
\\&=-2\f{H_{\xp}^{(2)}}{2\xp+1}\e -\f 2{2x+1}H_{\xp}^{(2)} \e -\f
1{2x+1}B_{p-2}(-x)\mod p.\endalign$$ Hence,
$$\align \sum_{n=0}^{p-1}V_n(x)&\e
p\sum_{k=0}^{p-1}\b xk\b{-1-x}k\f 1{2k+1}-p^3\sum_{k=0}^{p-1}\b
xk\b{-1-x}k\f{H_k^{(2)}}{2k+1}
\\&\e\f{1+2x'}{1+2x}p+\f{x'(x'+1)}{1+2x}p^3B_{p-2}(-x)+\f 1{1+2x}p^3B_{p-2}(-x)
\mod {p^4}.\endalign$$ By [S4, pp.216-217],
$$B_{p-2}\Ls 13\e 6U_{p-3}\mod p,\q B_{p-2}\Ls 16\e 30U_{p-3}\mod
p.$$ By [S2, Lemma 2.5], $B_{p-2}\ls 14\e 8E_{p-3}\mod p$. Now,
taking $x=-\f 13,-\f 14,-\f 16$ in the congruence for
$\sum_{n=0}^{p-1}V_n(x)\mod {p^4}$ and then applying the above
yields the congruences involving $G_n^{(3)},G_n^{(4)}$ and
$G_n^{(6)}$.
\par\q
\par{\bf Remark 4.2} Let $p$ be an odd prime. In [S20], the author proved
that
$$ \sum_{n=0}^{p-1}\f{V_n}{16^n}
\e 1+ \f 72p^3B_{p-3}\mod {p^4}.$$

\pro{Theorem 4.5} Let $p>3$ be a prime, $x\in\Bbb Z_p$, $x\not\e
0,-1\mod p$ and $x'=(x-\xp)/p$. Then
$$\sum_{n=0}^{p-1}(2n+1)V_n(x)\e  p^3\f{x'(x'+1)}{x(x+1)}
+2p^4\f{x'(x'+1)+1}{x(x+1)}H_{\xp}\mod {p^5}.
$$ Hence
$$\align
&\sum_{n=0}^{p-1}(2n+1)\f{V_n}{16^n} \e p^3+12p^4\qp 2\mod {p^5},
\\&\sum_{n=0}^{p-1}(2n+1)\f{V_n^{(3)}}{27^n}
\e p^3+\f{21}2p^4\qp 3\mod {p^5},
\\&\sum_{n=0}^{p-1}(2n+1)\f{V_n^{(4)}}{64^n} \e p^3+26p^4\qp
2\mod {p^5},
\\&\sum_{n=0}^{p-1}(2n+1)\f{V_n^{(6)}}{432^n}
\e p^3+\f{31}5p^4(4\qp 2+3\qp 3)\mod {p^5}\qtq{for}p>5.
\endalign$$
 \endpro
Proof. By (2.4),
$$\sum_{n=k}^{p-1}(2n+1)\b{n+k}{2k}=\f{p(p-k)}{k+1}\b{p+k}{2k}=\f{p^2}{k+1}
\cdot\f{(p^2-1^2)(p^2-2^2)\cdots(p^2-k^2)}{(2k)!}.$$ Thus,
$$\align\sum_{n=0}^{p-1}V_n(x)&=\sum_{n=0}^{p-1}(2n+1)
\sum_{k=0}^n\b nk\b{n+k}k(-1)^k\b xk\b{-1-x}k
\\&=\sum_{k=0}^{p-1}\b{x}k\b{-1-x}k(-1)^k\b{2k}k\sum_{n=k}^{p-1}
(2n+1)\b{n+k}{2k}
\\&=\sum_{k=0}^{p-1}\b{x}k\b{-1-x}k(-1)^k
\f{p^2}{k+1} \cdot\f{(p^2-1^2)(p^2-2^2)\cdots(p^2-k^2)}{k!^2}
\\&\e \sum_{k=0}^{p-1}\b{x}k\b{-1-x}k
\f{p^2}{k+1}\cdot\f{(-1)^k(-1^2)(-2^2)\cdots(-k^2)(1-p^2H_k^{(2)})}{k!^2}
\\&=p^2\sum_{k=0}^{p-1}\b{x}k\b{-1-x}k\f 1{k+1}-p^4
\sum_{k=0}^{p-1}\b{x}k\b{-1-x}k\f{H_k^{(2)}}{k+1}\mod {p^6}.
\endalign$$
Using the software Sigma we find the identity
$$\sum_{k=0}^n\b
nk\b{n+k}k(-1)^k\f{H_k^{(2)}}{k+1}=-\f{2H_n}{n(n+1)}.\tag 4.3$$
Thus,
$$\align \sum_{k=0}^{p-1}\b{x}k\b{-1-x}k\f{H_k^{(2)}}{k+1}
&=\sum_{k=0}^{p-1}\b{x}k\b{x+k}k(-1)^k\f{H_k^{(2)}}{k+1}
\\&\e\sum_{k=0}^{\xp}\b {\xp}k\b{\xp+k}k(-1)^k\f{H_k^{(2)}}{k+1}
\\&=-2\f{H_{\xp}}{\xp(\xp+1)}\e -2\f{H_{\xp}}{x(x+1)}\mod p.
\endalign$$
By [S8, Lemma 2.4] and [S12, Lemma 2.2],
$$\align\sum_{k=0}^{p-1}\b xk\b {-1-x}k\f p{k+1}
&= \b{x-1}{p-1}\b{-2-x}{p-1}=\f{p+x}{x+1}\b{x-1}{p-1}\b{-x-1}{p-1}
\\&\e\f{p+x}{x+1}\cdot p^2x'(x'+1)\Big(\f
1{\xp^2}-\f{p(1+2x')}{x^3}+\f{2pH_{\xp}}{x^2}\Big)
\\&\e \f{p^2x'(x'+1)}{x+1}\Big(\f
x{\xp^2}-2p\f{x'-xH_{\xp}}{x^2}\Big)
\\&=\f{p^2x'(x'+1)}{x+1}\Big(\f
x{(x-px')^2}-2p\f{x'}{x^2}+2p\f{H_{\xp}}x\Big)
\\&\e\f{p^2x'(x'+1)}{x+1}\Big(\f
1{x-2px'}-2p\f{x'}{x^2}+2p\f{H_{\xp}}x\Big)
\\&\e\f{p^2x'(x'+1)}{x+1}\Big(\f
{x+2px'}{x^2}-2p\f{x'}{x^2}+2p\f{H_{\xp}}x\Big)
\\&=p^2\f{x'(x'+1)}{x(x+1)}(1+2pH_{\xp})
 \mod {p^4}.
\endalign$$
Thus,
$$\align \sum_{n=0}^{p-1}(2n+1)V_n(x)
&\e p^3\f{x'(x'+1)}{x(x+1)}(1+2pH_{\xp})+2p^4\f{H_{\xp}}{x(x+1)}
\\&=p^3\f{x'(x'+1)}{x(x+1)}
+2p^4\f{x'(x'+1)+1}{x(x+1)}H_{\xp}\mod {p^5}.\endalign$$ By (2.17)
and (2.19)-(2.22),
$$\align &H_{\langle-\f 12\rangle_p}\e -2\qp 2\mod p,
\q H_{\langle-\f 13\rangle_p}\e -\f 32\qp 3\mod p,
\\&H_{\langle-\f 14\rangle_p}\e -3\qp 2\mod p,\q
H_{\langle-\f 16\rangle_p}\e -2\qp 2-\f 32\qp 3\mod p.
\endalign$$
Now taking $x=-\f 12,-\f 13,-\f 14,-\f 16$ in the above congruence
for $\sum_{n=0}^{p-1}(2n+1)V_n(x)\mod {p^5}$ and noting that
$x'(x'+1)=x(x+1)$ yields the remaining results.

\pro{Theorem 4.6} Let $p>3$ be a prime, $x\in\Bbb Z_p$  and
$x'=(x-\xp)/p$. Then
$$\sum_{n=0}^{p-1}(2n+1)(-1)^nV_n(x)\e  (-1)^{\xp}p+p^3(x'(x'+1)+1)
E_{p-3}(-x)\mod {p^4}.$$ Hence
$$\align
&\sum_{n=0}^{p-1}(2n+1)\f{V_n^{(3)}}{(-27)^n} \e (-1)^{[\f
p3]}p+7p^3U_{p-3}\mod {p^4},
\\&\sum_{n=0}^{p-1}(2n+1)\f{V_n^{(4)}}{(-64)^n} \e (-1)^{[\f p4]}
p+13p^3s_{p-3}\mod {p^4},
\\&\sum_{n=0}^{p-1}(2n+1)\f{V_n^{(6)}}{(-432)^n}
\e (-1)^{\f{p-1}2}p+\f{155}9p^3E_{p-3}\mod {p^4}.
\endalign$$
 \endpro
Proof. By (2.5),
$$\sum_{n=k}^{p-1}(2n+1)(-1)^n\b{n+k}{2k}=(p-k)\b{p+k}{2k}.$$
Thus,
$$\align &\b{2k}k\sum_{n=k}^{p-1}(2n+1)(-1)^n\b{n+k}{2k}
\\&=\f{p(p^2-1^2)(p^2-2^2)\cdots(p^2-k^2)}{k!^2} \e
(-1)^kp(1-p^2H_k^{(2)})\mod {p^5}.\endalign$$ It then follows that
$$\align\sum_{n=0}^{p-1}(2n+1)(-1)^nV_n(x)&=\sum_{n=0}^{p-1}(2n+1)(-1)^n
\sum_{k=0}^n\b nk\b{n+k}k(-1)^k\b xk\b{-1-x}k
\\&=\sum_{k=0}^{p-1}\b{x}k\b{-1-x}k(-1)^k\b{2k}k\sum_{n=k}^{p-1}
(2n+1)(-1)^n\b{n+k}{2k}
\\&\e\sum_{k=0}^{p-1}\b{x}k\b{-1-x}k(p-p^3H_k^{(2)})\mod {p^5}.
\endalign$$
By [T2],
$$\sum_{n=1}^k\b
nk\b{n+k}k(-1)^kH_k^{(2)}=2(-1)^{n+1}\sum_{k=1}^n\f{(-1)^k}{k^2}.\tag
4.4$$ Thus, appealing to (2.12),
$$\align &\sum_{k=0}^{p-1}\b{x}k\b{-1-x}kH_k^{(2)}
\\&=\sum_{k=0}^{p-1}\b{x}k\b{x+k}k(-1)^kH_k^{(2)} \e \sum_{k=0}^{\xp}
\b{\xp}k\b{\xp+k}k(-1)^kH_k^{(2)} \\&=2(-1)^{\xp+1} \sum_{k=1}^{\xp}
\f{(-1)^k}{k^2}\e -E_{p-3}(-x)\mod p.\endalign$$ From the above and
(2.13),
$$\align \sum_{n=0}^{p-1}(2n+1)(-1)^nV_n(x)&\e
p\sum_{k=0}^{p-1}\b
xk\b{-1-x}k-p^3\sum_{k=0}^{p-1}\b{x}k\b{-1-x}kH_k^{(2)}
\\&\e (-1)^{\xp}p+p^3x'(x'+1)E_{p-3}(-x)+p^3E_{p-3}(-x)\mod {p^4}.
\endalign$$
Now, taking $x=-\f 13,-\f 14,\f 16$ and then applying (2.14)-(2.15)
yields the congruences involving $V_n^{(3)},V_n^{(4)}$ and
$V_n^{(6)}$. The proof is now complete.
\par\q
\par{\bf Remark 4.3} Let $p>3$ be a prime.
 In [Su3], Z.W. Sun conjectured that
$$\align
&\sum_{n=0}^{p-1}(n+1)\f{V_n}{8^n}\e (-1)^{\f{p-1}2}p+5p^3
E_{p-3}\mod {p^4},\tag 4.5
\\&\sum_{n=0}^{p-1}(2n+1)\f{V_n}{(-16)^n}\e (-1)^{\f{p-1}2}p+3p^3
E_{p-3}\mod {p^4}.\tag 4.6\endalign$$ These congruences were
recently proved by Wang [W]. (4.6) can be deduced from Theorem 4.6
(with $x=-\f 12$).

\pro{Theorem 4.7} Let $p$ be an odd prime, $x\in\Bbb Z_p$ and
$x\not\e 0\mod p$. Then
$$\sum_{n=0}^{p-1}(-1)^nV_n(x)\e
\cases \sum_{k=0}^{(p-1)/2}\b{2k}k\b xk\b{-1-x}k\f 1{4^k}\mod {p^2}
&\t{if $\xp$ is even,}\\0\mod p &\t{if $\xp$ is odd.}
\endcases$$
Hence, for $p>3$,
$$\align &\sum_{n=0}^{p-1}\f{V_n^{(3)}}{(-27)^n}\e
\cases 4x^2-2p\mod {p^2}&\t{if $p=x^2+3y^2\e 1\mod 3$,}
\\0\mod p&\t{if $p\e 2\mod 3$,}\endcases
\\&\sum_{n=0}^{p-1}\f{V_n^{(4)}}{(-64)^n}\e
\cases 4x^2-2p\mod {p^2}&\t{if $p=x^2+2y^2\e 1,3\mod 8$,}
\\0\mod p&\t{if $p\e 5,7\mod 8$,}\endcases
\\&\sum_{n=0}^{p-1}\f{V_n^{(6)}}{(-432)^n}\e \cases
 (-1)^{[\f p3]}(4x^2-2p)\mod
{p^2}&\t{if $p=x^2+4y^2\e 1\mod 4$,}
\\0\mod p&\t{if $p\e 3\mod 4$.}
\endcases\endalign$$
\endpro
Proof. By Theorem 4.1, $(-1)^nV_n(x)=\sum_{k=0}^n\b
nk\b{n+k}k(-1)^kG_k(x)$. Thus, applying Lemma 2.4 gives
$$\sum_{n=0}^{p-1}(-1)^nV_n(x)\e\sum_{k=0}^{p-1}
\f p{2k+1}G_k(x)=G_{\f{p-1}2}(x)+p\sum_{k=0}^{(p-3)/2}
\Big(\f{G_k(x)}{2k+1}+\f{G_{p-1-k}(x)}{2(p-1-k)+1}\Big)\mod {p^3}.$$
By [S18, Theorem 3.1] and Theorem 3.4, for $k=0,1,\ldots,p-1$,
$$G_k(x)\e G_{p-1}(x)G_{p-1-k}(x)\e (-1)^{\xp}G_{p-1-k}(x)\mod p.
\tag 4.7$$ Thus,
$$\align \sum_{n=0}^{p-1}(-1)^nV_n(x)
&\e G_{\f{p-1}2}(x)+p\sum_{k=0}^{(p-3)/2}
\Big(\f{G_k(x)}{2k+1}-\f{G_{p-1-k}(x)}{2k+1}\Big) \\& \e
G_{\f{p-1}2}(x)+p\big(1-(-1)^{\xp}\big)\sum_{k=0}^{(p-3)/2}\f
{G_k(x)}{2k+1}\mod {p^2}.\endalign$$ Now applying Theorem 3.6 yields
the result.

\pro{Theorem 4.8} Let $p>3$ be a prime,  $m,x\in\Bbb Z_p$ and
$m\not\e 0\mod p$. Then
$$\sum_{n=0}^{p-1}\f{V_n(x)}{m^n}
\e\Big(\sum_{k=0}^{p-1}\f{G_k(x)}{m^k}\Big)^2\e
\sum_{k=0}^{p-1}\f{\b{2k}k\b xk \b{-1-x}k}{(2-m-1/m)^k}\mod
p\qtq{for}m\not\e \pm 1\mod p.$$ Hence
$$\align &\sum_{n=0}^{p-1}\f{V_n^{(3)}}{m^n}
\e
\Big(\sum_{k=0}^{p-1}\f{G_k^{(3)}}{m^k}\Big)^2\e\sum_{k=0}^{p-1}\f{\b{2k}k^2\b{3k}k}{(54-m-729/m)^k}\mod
p \ \t{for}\ m\not\e \pm 27\mod p,
\\&\sum_{n=0}^{p-1}\f{V_n^{(4)}}{m^n}\e
\Big(\sum_{k=0}^{p-1}\f{G_k^{(4)}}{m^k}\Big)^2
\e\sum_{k=0}^{p-1}\f{\b{2k}k^2\b{4k}{2k}}{(128-m-4096/m)^k}\mod p \
\t{for}\ m\not\e \pm 64\mod p,
\\&\sum_{n=0}^{p-1}\f{V_n^{(6)}}{m^n}\e
\Big(\sum_{k=0}^{p-1}\f{G_k^{(6)}}{m^k}\Big)^2
\e\sum_{k=0}^{p-1}\f{\b{2k}k\b{3k}k\b{6k}{3k}}{(864-m-186624/m)^k}\mod
p\ \t{for}\ m\not\e \pm 432\mod p.
\endalign$$
\endpro
Proof. By Theorem 4.1, $(-1)^nV_n(x)=\sum_{k=0}^n\b
nk\b{n+k}k(-1)^kG_k(x)$. Thus, taking $c_k=(-1)^kG_k(x)$ in Lemma
2.3 gives
$$\sum_{n=0}^{p-1}V_n(x)(-u)^n\e
\sum_{k=0}^{p-1}\b{2k}k\Ls {-u}{(1-u)^2}^kG_k(x)\mod p\qtq{for}
u\not\e 1\mod p.$$ Replacing $u$ with $-1/m$ and then applying
Theorem 3.7 yields
$$\align \sum_{n=0}^{p-1}\f{V_n(x)}{m^n}
&\e\sum_{k=0}^{p-1}\f{\b{2k}kG_k(x)}{(m+1/m+2)^k} \e\Ls{(m+\f
1m+2)(m+\f 1m-2)}p \sum_{k=0}^{p-1}\f{\b{2k}k\b xk\b{-1-x}k}
{(2-m-1/m)^k}
\\&= \sum_{k=0}^{p-1}\f{\b{2k}k\b xk\b{-1-x}k}
{(2-m-1/m)^k} \mod p\qtq{for}m\not\e \pm 1 \mod p.\endalign$$ From
Theorem 3.7, [S11, Theorem 2.2] and the above,
$$ \Big(\sum_{k=0}^{p-1}\f{G_k(x)}{m^k}\Big)^2\e
\Big(\sum_{k=0}^{p-1} \f{\b xk\b{-1-x}k}{(1-m)^k}\Big)^2 \e
\sum_{k=0}^{p-1}\f{\b{2k}k\b xk\b{-1-x}k} {(2-m-1/m)^k} \e
\sum_{n=0}^{p-1}\f{V_n(x)}{m^n}\mod p.$$ Now taking $x=-\f 13,-\f
14,-\f 16$ and replacing $m$ with $\f m{27},\f{m}{64},\f m{432}$
yields the remaining results.
\par\q
\par Now we pose some challenging conjectures involving $V_n,V_n^{(3)},V_n^{(4)}$
and $V_n^{(6)}$.
 \pro{Conjecture 4.1}
Let $p>3$ be a prime. Then
$$\align &
V_p\e 8+40p^3B_{p-3}\mod {p^4},\ V_p^{(3)}\e 15+132p^3B_{p-3}\mod
{p^4},
\\&V_p^{(4)}\e 40+704p^3B_{p-3}\mod {p^4},\ V_p^{(6)}\e 312+
16120p^3B_{p-3}\mod {p^4}, \\&V_{2p}\e V_2+2624p^3B_{p-3}\mod {p^4},
\q V_{2p}^{(3)}\e V_2^{(3)}+16416p^3B_{p-3}\mod {p^4},
\\&V_{2p}^{(4)}\e V_2^{(4)}+233984p^3B_{p-3}\mod {p^4}.
\endalign$$
\endpro
\pro{Conjecture 4.2} Let $p>3$ be a prime. Then
$$\align & V_{p-1}^{(3)}\e 729^{p-1}-\f{92}{27}p^3B_{p-3}\mod {p^4},
\\& V_{p-1}^{(4)}\e 4096^{p-1}-\f{17}2p^3B_{p-3}\mod {p^4}, \q
\\&V_{p-1}^{(6)}\e 186624^{p-1}-\f{1705}{54}p^3B_{p-3} \mod {p^4},
\\&V_{2p-1}^{(3)}\e 27^{4(p-1)}V_1^{(3)}-\f{1504}3p^3B_{p-3}\mod
{p^4}, \\&V_{2p-1}^{(4)}\e 64^{4(p-1)}V_1^{(4)}-3196p^3B_{p-3}\mod
{p^4},
\\&V_{2p-1}^{(6)}\e 432^{4(p-1)}V_1^{(6)}-\f {264940}3p^3B_{p-3}\mod
{p^4}.
\endalign$$
\endpro

\pro{Conjecture 4.3} Let $p$ be a prime with $p>3$. Then
$$\sum_{n=0}^{p-1}\f{V_n^{(3)}}{(-27)^n}
\e \cases 4x^2-2p-\f{p^2}{4x^2}\mod {p^3}&\t{if $p=x^2+3y^2\e 1\mod
3$,}
\\\f 74p^2\b{\f{p-1}2}{\f{p-5}6}^{-2}\mod {p^3}
&\t{if $p\e 2\mod 3$}.\endcases$$
\endpro

\pro{Conjecture 4.4} Let $p$ be a prime with $p>3$. Then
$$\align&\sum_{n=0}^{p-1}\f{V_n^{(3)}}{3^n}
\e \cases L^2-2p-\f{p^2}{L^2}\mod {p^3}&\t{if $p\e 1\mod 3$ and so
$4p=L^2+27M^2$,}
\\\f {15}2p^2\b{[\f{2p}3}{[\f p3]}^{-2}\mod {p^3}
&\t{if $p\e 2\mod 3$},\endcases
\\&\sum_{n=0}^{p-1}\f{V_n^{(3)}}{243^n}
\e \cases L^2-2p-\f{p^2}{L^2}\mod {p^3}&\t{if $p\e 1\mod 3$ and so
$4p=L^2+27M^2$,}
\\-\f 12p^2\b{[\f{2p}3}{[\f p3]}^{-2}\mod {p^3}
&\t{if $p\e 2\mod 3$}.\endcases
\endalign$$
\endpro

\pro{Conjecture 4.5} Let $p$ be a prime with $p\not=2,3,7$. Then
$$\align&\sum_{n=0}^{p-1}V_n^{(4)}
\e \cases 4x^2-2p-\f{p^2}{4x^2}\mod {p^3}&\t{if $p=x^2+7y^2\e
1,2,4\mod 7$,}
\\\f{149}{12}p^2\b{3[p/7]}{[p/7]}^{-2}\mod {p^3}&\t{if $p\e 3\mod
7$,}
\\\f{447}{64}p^2\b{3[p/7]}{[p/7]}^{-2}\mod {p^3}&\t{if $p\e 5\mod
7$,}
\\\f{3725}{5808}p^2\b{3[p/7]}{[p/7]}^{-2}\mod {p^3}&\t{if $p\e 6\mod
7$,}\endcases
\\&\sum_{n=0}^{p-1}\f{V_n^{(4)}}{4096^n}
\e \cases 4x^2-2p-\f{p^2}{4x^2}\mod {p^3}&\t{if $p=x^2+7y^2\e
1,2,4\mod 7$,}
\\-\f{5}{96}p^2\b{3[p/7]}{[p/7]}^{-2}\mod {p^3}&\t{if $p\e 3\mod
7$,}
\\-\f{15}{512}p^2\b{3[p/7]}{[p/7]}^{-2}\mod {p^3}&\t{if $p\e 5\mod
7$,}
\\-\f{125}{46464}p^2\b{3[p/7]}{[p/7]}^{-2}\mod {p^3}&\t{if $p\e 6\mod
7$.}\endcases
\endalign$$
\endpro

\pro{Conjecture 4.6} Let $p$ be a prime with $p>3$. Then
$$\align &\sum_{n=0}^{p-1}\f{V_n^{(4)}}{16^n}
\e \cases 4x^2-2p-\f{p^2}{4x^2}\mod {p^3}&\t{if $p=x^2+3y^2\e 1\mod
3$,}
\\5p^2\b{\f{p-1}2}{\f{p-5}6}^{-2}\mod {p^3}
&\t{if $p\e 2\mod 3$},\endcases
\\&\sum_{n=0}^{p-1}\f{V_n^{(4)}}{256^n}
\e \cases 4x^2-2p-\f{p^2}{4x^2}\mod {p^3}&\t{if $p=x^2+3y^2\e 1\mod
3$,}
\\-p^2\b{\f{p-1}2}{\f{p-5}6}^{-2}\mod {p^3}
&\t{if $p\e 2\mod 3$}.\endcases
\endalign$$
\endpro

\pro{Conjecture 4.7} Let $p$ be a prime with $p>3$. Then
$$\align &\sum_{n=0}^{p-1}\f{V_n^{(4)}}{(-8)^n}
\e \cases 4x^2-2p-\f{p^2}{4x^2}\mod {p^3}&\t{if $p=x^2+4y^2\e 1\mod
4$,}
\\\f{41}{12}p^2\b{\f{p-3}2}{\f{p-3}4}^{-2}\mod {p^3}
&\t{if $p\e 3\mod 4$},\endcases
\\&\sum_{n=0}^{p-1}\f{V_n^{(4)}}{(-512)^n}
\e \cases 4x^2-2p-\f{p^2}{4x^2}\mod {p^3}&\t{if $p=x^2+4y^2\e 1\mod
4$,}
\\\f 5{12}p^2\b{\f{p-3}2}{\f{p-3}4}^{-2}\mod {p^3}
&\t{if $p\e 3\mod 4$}.\endcases
\endalign$$
\endpro

\pro{Conjecture 4.8} Let $p$ be a prime with $p>3$. Then
$$(-1)^{\f{p-1}2}\sum_{n=0}^{p-1}\f{V_n^{(4)}}{(-64)^n}
\e \cases 4x^2-2p-\f{p^2}{4x^2}\mod {p^3}&\t{if $p=x^2+2y^2\e
1,3\mod 8$,}
\\-\f{13}9p^2\b{[\f p4]}{[\f p8]}^{-2}\mod {p^3}&\t{if $p\e 5\mod
8$,}
\\-\f{13}2p^2\b{[\f p4]}{[\f p8]}^{-2}\mod {p^3}&\t{if $p\e 7\mod
8$.}\endcases$$
\endpro

\pro{Conjecture 4.9} Let $p$ be a prime with $p>3$. Then
$$(-1)^{[\f p3]}\sum_{n=0}^{p-1}\f{V_n^{(6)}}{(-432)^n}
\e \cases 4x^2-2p-\f{p^2}{4x^2}\mod {p^3}&\t{if $p=x^2+4y^2\e 1\mod
4$,}
\\-\f{31}{12}p^2\b{\f{p-3}2}{\f{p-3}4}^{-2}\mod {p^3}
&\t{if $p\e 3\mod 4$}.\endcases
$$\endpro

\pro{Conjecture 4.10} Let $p>3$ be a prime and $m,r\in\Bbb Z^+$.
Then
$$V_{mp^r}^{(3)}\e V_{mp^{r-1}}^{(3)}\pmod {p^{3r}},
\ V_{mp^r}^{(4)}\e V_{mp^{r-1}}^{(4)}\pmod {p^{3r}},\
V_{mp^r}^{(6)}\e V_{mp^{r-1}}^{(6)}\pmod {p^{3r}}.$$ Moreover,
$$\align &\f{V_{mp^r}^{(3)}-V_{mp^{r-1}}^{(3)}}{p^{3r}}
\e \f{V_{mp}^{(3)}-V_m^{(3)}}{p^3}\mod p,\q
\f{V_{mp^r}^{(4)}-V_{mp^{r-1}}^{(4)}}{p^{3r}} \e
\f{V_{mp}^{(4)}-V_m^{(4)}}{p^3}\mod p,
\\&\f{V_{mp^r}^{(6)}-V_{mp^{r-1}}^{(6)}}{p^{3r}}
\e \f{V_{mp}^{(6)}-V_m^{(6)}}{p^3}\mod p
 .\endalign$$
\endpro

\pro{Conjecture 4.11} Let $p>3$ be a prime and $m,r\in\Bbb Z^+$.
Then
$$ \align &V_{mp^r-1}^{(3)}\e 729^{mp^{r-1}(p-1)}V_{mp^{r-1}-1}^{(3)}
\mod {p^{3r}},
\\&V_{mp^r-1}^{(4)}\e 4096^{mp^{r-1}(p-1)}
V_{mp^{r-1}-1}^{(4)} \mod {p^{3r}},
\\&V_{mp^r-1}^{(6)}\e 186624^{mp^{r-1}(p-1)}
V_{mp^{r-1}-1}^{(6)} \mod {p^{3r}}  . \endalign$$ Moreover,
$$ \align &\f{V_{mp^r-1}^{(3)}-729^{mp^{r-1}(p-1)}
V_{mp^{r-1}-1}^{(3)}}{p^{3r}} \e
\f{V_{mp-1}^{(3)}-729^{m(p-1)}V_{m-1}^{(3)}}{p^3}\mod p,
\\&\f{V_{mp^r-1}^{(4)}-4096^{mp^{r-1}(p-1)}
V_{mp^{r-1}-1}^{(4)}}{p^{3r}} \e
\f{V_{mp-1}^{(4)}-4096^{m(p-1)}V_{m-1}^{(4)}}{p^3}\mod p,
\\&\f{V_{mp^r-1}^{(6)}-186624^{mp^{r-1}(p-1)}
V_{mp^{r-1}-1}^{(6)}}{p^{3r}} \e
\f{V_{mp-1}^{(6)}-186624^{m(p-1)}V_{m-1}^{(6)}}{p^3}\mod
p.\endalign$$
\endpro
 \pro{Conjecture 4.12} Let $x$ be a real number and $n\in\Bbb Z^+$.
 If $-1<x<0$, then
$V_n(x)^2<V_{n+1}(x)V_{n-1}(x)$. If $x\notin [-1,0]$, then
$V_n(x)^2>V_{n+1}(x)V_{n-1}(x)$.
\endpro

\section*{5. Congruences involving other Ap\'ery-like sequences}
\par In this section, we use Lemmas 2.1 and 2.2 to deduce
congruences involving $a_n,f_n$ and $Q_n$, and pose many conjectures
on congruences involving known Ap\'ery-like sequences.

\pro{Lemma 5.1} Let $n$ be a nonnegative integer. Then
$$\align &f_n=\sum_{k=0}^n\b nk(-1)^{n-k}a_k=\sum_{k=0}^n\b nk8^{n-k}Q_k,
\\&Q_n=\sum_{k=0}^n\b nk(-9)^{n-k}a_k,
\\&a_n=\sum_{k=0}^n\b nkf_k=\sum_{k=0}^n\b nk9^{n-k}Q_k.\endalign$$
\endpro
Proof. By [St, (38)], $a_n=\sum_{k=0}^n\b nkf_k$.  Applying the
binomial inversion formula gives $f_n=\sum_{k=0}^n\b
nk(-1)^{n-k}a_k$. Since $\f{Q_n}{(-8)^n}=\sum_{k=0}^n\b
nk(-1)^k\f{f_k}{8^k}$, applying the binomial inversion formula gives
$\f{f_n}{8^n}=\sum_{k=0}^n\b nk\f{Q_k}{8^k}$. Also,
$$\align &\sum_{k=0}^n\b nk(-9)^{n-k}a_k
\\&=\sum_{k=0}^n\b nk(-1)^{n-k}a_k\sum_{r=0}^{n-k}\b{n-k}r8^r
=\sum_{r=0}^n\b nr(-8)^r\sum_{k=0}^{n-r}\b{n-r}k(-1)^{n-r-k}a_k
\\&=\sum_{r=0}^n\b nr(-8)^rf_{n-r}=Q_n.\endalign$$
Applying the binomial inversion formula yields $a_n=\sum_{k=0}^n\b
nk9^{n-k}Q_k$, which completes the proof.
 \pro{Theorem 5.1} Suppose
that $p$ is an odd prime, $m\in\Bbb Z_p$ and $m\not\e 0,1\mod p$.
Then
$$\align\sum_{k=0}^{p-1}\f{a_k}{m^k}
\e \cases \sum_{k=0}^{p-1}\f{\b{2k}k\b{3k}k}{((m-3)^3/(m-1))^k}\mod
p&\t{if $m\not\e 3\mod p$,}
\\\sum_{k=0}^{p-1}\f{\b{2k}k\b{3k}k}{((m+3)^3/(m-1)^2)^k}\mod
p&\t{if $m\not\e -3\mod p$}\endcases
\endalign$$
and
$$\sum_{k=0}^{p-1}\f{Q_k}{(-8m)^k}
\e \cases \sum_{k=0}^{p-1}\f{\b{2k}k\b{3k}k}
{((4m-3)^3/(m-1))^k}\mod p&\t{if $m\not\e \f 34\mod p$,}
\\\sum_{k=0}^{p-1}\f{\b{2k}k\b{3k}k}
{(-(2m-3)^3/(m-1)^2)^k}\mod p&\t{if $m\not\e \f 32\mod p$.}
\endcases$$
\endpro
Proof. By Lemma 5.1, $a_n=\sum_{k=0}^n\b nkf_k$. From Lemma 2.1,
$\sum_{k=0}^{p-1}\f{a_k}{m^k}\e\sum_{k=0}^{p-1}\f{f_k}{(m-1)^k}$
$\mod p$. Now applying [S16, Theorem 2.12 and Lemmas 2.4 (with $z=\f
1{m-1}$)] yields the first result.
 Since $\f{Q_n}{(-8)^n}=\sum_{k=0}^n\b nk\f{f_k}{(-8)^k}$,
applying Lemma 2.1 gives $\sum_{k=0}^{p-1}\f{Q_k}{(-8m)^k} \e
\sumkp\f{f_k}{(-8(m-1))^k}\mod p$. From [S16,  Theorem 2.12 and
Lemma 2.4 (with $z=\f 1{-8(m-1)}$)] we deduce the remaining part.
\par\q
\pro{Theorem 5.2} Suppose that $p$ is a prime with $p>3$. Then
$$\sum_{n=0}^{p-1}\f{Q_n}{(-6)^n}\e \sum_{n=0}^{p-1}\f{Q_n}{(-12)^n} \e\cases 2x\mod
p&\t{if $3\mid p-1$, $p=x^2+3y^2$ and $3\mid x-1$,}
\\0\mod p&\t{if $p\e 2\mod 3$.}
\endcases$$
\endpro
Proof. Putting $m=\f 34$ and $\f 32$ in Theorem 5.1 yields
$$\sum_{n=0}^{p-1}\f{Q_n}{(-6)^n}\e \f{Q_n}{(-12)^n}\e
\sum_{k=0}^{p-1}\f{\b{2k}k\b{3k}k}{54^k}\mod p.$$ Now applying [S8,
Theorem 3.4] yields the result.
 \pro{Theorem 5.3} Let $p$ be an odd
prime, $m\in\Bbb Z_p$ and $(m+2)(m-2)\not\e 0\mod p$. Then
$$\align &\sum_{k=0}^{p-1}\b{2k}k\f{a_k}{(m+2)^k}
\e\Ls{(m+2)(m-2)}p\sum_{k=0}^{p-1}\b{2k}k\f{f_k}{(m-2)^k}\mod p,\tag
5.1
\\&\sum_{k=0}^{p-1}\b{2k}k\f{Q_k}{(-8(m+2))^k}
\e\Ls{(m+2)(m-2)}p\sum_{k=0}^{p-1}\b{2k}k\f{f_k}{(-8(m-2))^k}\mod p,
\tag 5.2
\\&\sum_{k=0}^{p-1}\b{2k}k\f{Q_k}{(-9(m+2))^k}
\e\Ls{(m+2)(m-2)}p\sum_{k=0}^{p-1}\b{2k}k\f{a_k}{(-9(m-2))^k}\mod p
\tag 5.3
\endalign$$
and for $9m-14\not\e 0\mod p$,
$$\sum_{k=0}^{p-1}\b{2k}k\f{Q_k}{(-9(m+2))^k}
\e\Ls{(m+2)(9m-14)}p \sum_{k=0}^{p-1}\b{2k}k\f{f_k}{(-9m+14)^k}\mod
p.\tag 5.4$$
\endpro
Proof. We first note that $p\mid \b{2k}k$ for $\f {p+1}2\le k\le
p-1$. By Lemma 5.1, taking $u_k=f_k$ and $v_k=a_k$ in Lemma 2.2
gives (5.1). Since $\f{Q_n}{(-8)^n}=\sum_{k=0}^n \b
nk\f{f_k}{(-8)^k}$, taking $u_k=\f{f_k}{(-8)^k}$ and
$v_k=\f{Q_k}{(-8)^k}$ in Lemma 2.2 gives (5.2). By Lemma 5.1, taking
$u_k=\f{a_k}{(-9)^k}$ and $v_k=\f{Q_k}{(-9)^k}$ in Lemma 2.2 yields
(5.3). Combining (5.3) with (5.1) yields (5.4).
\par\q
\pro{Theorem 5.4} Suppose that $p$ is a prime with $p>5$. Then
$$\align (-1)^{\f{p-1}2}\sum_{n=0}^{p-1}\b{2n}n\f{a_n}{54^n}&
\e \sum_{n=0}^{p-1}\b{2n}n\f{Q_n}{18^n}\e
\sum_{n=0}^{p-1}\b{2n}n\f{Q_n}{(-36)^n}
\\&\e\cases 4x^2\mod p&\t{if $p\e 1\mod 3$ and so $p=x^2+3y^2$,}
\\0\mod p&\t{if $p\e 2\mod 3$.}\endcases\endalign$$
\endpro
Proof. Taking $m=52$ in (5.1) and then applying [S16, Theorem 2.2]
gives
$$\align \sum_{n=0}^{p-1}\b{2n}n\f{a_n}{54^n}&\e \Ls 3p
\sum_{k=0}^{p-1}\b{2k}k\f{f_k}{50^k}\\&\e \cases
(-1)^{\f{p-1}2}4x^2\mod p&\t{if $3\mid p-1$ and so $p=x^2+3y^2$,}
\\0\mod p&\t{if $p\e 2\mod 3$.}\endcases\endalign$$
Taking $m=-4$ in (5.3) gives $\sum_{n=0}^{p-1}\b{2n}n\f{Q_n}{18^n}
\e \ls 3p\sum_{n=0}^{p-1}\b{2n}n\f{a_n}{54^n}\mod p.$ Taking $m=2$
in (5.4) and applying the congruence for
$\sum_{k=0}^{p-1}\b{2k}k\f{f_k}{(-4)^k}\mod p$ (see [S16, p.124])
yields
$$\align \sum_{n=0}^{p-1}\b{2n}n\f{Q_n}{(-36)^n}&\e
\sum_{k=0}^{p-1}\b{2k}k\f{f_k}{(-4)^k}\\&\e \cases 4x^2\mod p&\t{if
$3\mid p-1$ and so $p=x^2+3y^2$,}
\\0\mod p&\t{if $p\e 2\mod 3$.}\endcases\endalign$$
This proves the theorem.

\pro{Theorem 5.5} Suppose that $p$ is a prime such that $p\e
1,19\mod{30}$ and so $p=x^2+15y^2$. Then
$$\sum_{n=0}^{p-1}\b{2n}n\f{a_n}{9^n}\e
\sum_{n=0}^{p-1}\b{2n}n\f{a_n}{(-45)^n}\e
\sum_{n=0}^{p-1}\b{2n}n\f{Q_n}{(-27)^n}\e
\sum_{n=0}^{p-1}\b{2n}n\f{Q_n}{(-81)^n}\e 4x^2\mod p.$$
\endpro
Proof. Taking $m=7$ in (5.1) and then applying [S16, Theorem 2.5]
gives
$$\sum_{n=0}^{p-1}\b{2n}n\f{a_n}{9^n}\e \Ls 5p
\sum_{k=0}^{p-1}\b{2k}k\f{f_k}{5^k}\e 4x^2\mod p.$$ Taking $m=-47$
in (5.1) and then applying [S16, Theorem 2.4] gives
$$\sum_{n=0}^{p-1}\b{2n}n\f{a_n}{(-45)^n}\e \Ls 5p
\sum_{k=0}^{p-1}\b{2k}k\f{f_k}{(-49)^k}\e 4x^2\mod p.$$ Taking
$m=1,7$ in (5.3) gives
$$\align&\sum_{n=0}^{p-1}\b{2n}n\f{Q_n}{(-27)^n}\e
\sum_{k=0}^{p-1}\b{2k}k\f{a_k}{9^k}\mod p,
\\&\sum_{n=0}^{p-1}\b{2n}n\f{Q_n}{(-81)^n}\e
\sum_{k=0}^{p-1}\b{2k}k\f{a_k}{(-45)^k}\mod p. \endalign$$ Now
combining the above proves the theorem.

\pro{Theorem 5.6} Suppose that $p$ is a prime with $p>3$. Then
$$\align&\sum_{n=0}^{p-1}\b{2n}n\f{Q_n}{(-32)^n}
\e \sum_{n=0}^{p-1}\b{2n}n\f{Q_n}{64^n} \e \cases 4x^2\mod p&\t{if
$p=x^2+2y^2\e 1,3\mod 8$,}\\0\mod p&\t{if $p\e 5,7\mod 8$,}
\endcases
\\&\sum_{n=0}^{p-1}\b{2n}n\f{a_n}{20^n}
\e\sum_{n=0}^{p-1}\b{2n}n\f{Q_n}{(-16)^n}\e 4x^2\mod
p\qtq{for}p=x^2+5y^2\e 1,9\mod{20},
\\&\sum_{n=0}^{p-1}\b{2n}n\f{Q_n}{(-48)^n}\e
\cases 4x^2\mod p&\t{if $p=x^2+9y^2\e 1\mod{12}$,}
\\0\mod p&\t{if $p\e 11\mod{12}$.}
\endcases
\endalign$$
\endpro
Proof. Taking $m=\f{14}9, -\f{82}9$ in (5.3) yields
$$\align &\sum_{n=0}^{p-1}\b{2n}n\f{Q_n}{(-32)^n}\e
\sum_{n=0}^{p-1}\b{2n}n\f{a_n}{4^n}\mod p,
\\&\sum_{n=0}^{p-1}\b{2n}n\f{Q_n}{64^n}\e
\sum_{n=0}^{p-1}\b{2n}n\f{a_n}{100^n}\mod p.\endalign$$ Now applying
[S16, Theorem 2.14] and [S9, Theorem 5.6] yields the first
congruence. Taking $m=0$ in (5.2), $m=18$ in (5.1) and then applying
[S16, Theorem 2.10] yields the second congruence. Taking $m=\f{10}3$
in (5.3) and then applying [S9, Theorem 4.3] gives the third
congruence.
\par\q
\par Based on calculations by Maple, we pose the following
conjectures.
\pro{Conjecture 5.1} Let $p$ be a prime with $p>3$.
Then
$$\align&A_{p}\e A_1-\f{14}3p^3B_{p-3}\mod {p^4},
\q A_{2p}\e A_2-\f{1648}3p^3B_{p-3}\mod {p^4}, \q
 \\&A_{3p}\e
A_3-36738p^3B_{p-3}\mod {p^4},\q A_p'\e A_1'-\f 53p^3B_{p-3}\mod
{p^4},
\\&A_{2p}'\e A_2'-\f{280}3p^3B_{p-3}\mod {p^4},\q
A_{3p}'\e A_3'-2475p^3B_{p-3}\mod {p^4},
\\&b_p\e b_1-6p^3B_{p-3}\mod {p^4},\qq b_{2p}\e b_2+144p^3B_{p-3}\mod {p^4},
\\&b_{3p}\e b_3-1566p^3B_{p-3}\mod {p^4}, \q T_p\e T_1-p^3B_{p-3}\mod{p^4},
\\&T_{2p}\e T_2-136p^3B_{p-3}\mod {p^4},\q
T_{3p}\e T_3-6696p^3B_{p-3}\mod {p^4},
\\&D_p\e D_1+\f{16}3p^3B_{p-3}\mod {p^4},\q
D_{2p}\e D_2+\f{448}3p^3B_{p-3}\mod {p^4},
\\&D_{3p}\e D_3+3168p^3B_{p-3}\mod {p^4},\q
f_p\e f_1+\f 12p^3B_{p-3}\mod {p^4},
\\&f_{2p}\e f_2-8p^3B_{p-3}\mod {p^4},\q
f_{3p}\e f_3-189p^3B_{p-3}\mod {p^4}.
\endalign$$
\endpro
\par We remark that Z.W. Sun informed the author he
conjectured the congruence $A_{p}\e A_1-\f{14}3p^3B_{p-3}\mod {p^4}$
earlier.

\pro{Conjecture 5.2} Let $p$ be a prime with $p>3$. Then
$$\align &A_{2p-1}\e A_1+\f{16}3p^3B_{p-3}\mod {p^4},
\\&A_{2p-1}'\e A_1'+\f{200}3p^3B_{p-3}\mod {p^4},
\\&f_{2p-1}\e 8^{2(p-1)}f_1+17p^3B_{p-3}\mod {p^4},
\\&T_{2p-1}\e 16^{2(p-1)}T_1-6p^3B_{p-3}\mod {p^4},
\\&D_{2p-1}\e 64^{2(p-1)}D_1-\f{44}3p^3B_{p-3}\mod {p^4},
\\&b_{2p-1}\e 81^{2(p-1)}b_1+\f{16}9p^3B_{p-3}\mod {p^4}.
\endalign$$
\endpro
\pro{Conjecture 5.3} Let $p$ be a prime with $p>3$. Then
$$\align
&a_p\e a_1+3p^2\Ls p3U_{p-3}\mod {p^3},\q a_{2p}\e a_2+36p^2\Ls
p3U_{p-3}\mod {p^3},
\\&a_{3p}\e a_3+405p^2\Ls
p3U_{p-3}\mod {p^3},\q W_p\e W_1-9p^2\Ls p3U_{p-3}\mod{p^3},
\\&W_{2p}\e W_2+108p^2\Ls p3U_{p-3}\mod {p^3},\q
W_{3p}\e W_3-729p^2\Ls p3U_{p-3}\mod {p^3},
\\&Q_p\e Q_1-30p^2\Ls p3U_{p-3}\mod {p^3},\q
Q_{2p}\e Q_2+720p^2\Ls p3U_{p-3}\mod {p^3},
\\&Q_{3p}\e Q_3-11340p^2\Ls p3U_{p-3}\mod {p^3}.
\endalign$$
\endpro

 \pro{Conjecture 5.4} Let $p>3$ be a prime. Then
$$\align &a_{2p-1}\e (-1)^{[\f p3]}9^{2(p-1)}a_1+20p^2U_{p-3}\mod
{p^3}, \\&W_{2p-1}\e (-1)^{[\f p3]}27^{2(p-1)}W_1-12p^2U_{p-3}\mod
{p^3},
\\&Q_{2p-1}\e (-1)^{[\f p3]}72^{2(p-1)}Q_1-70p^2U_{p-3}\mod
{p^3}.
\endalign$$
\endpro

 \pro{Conjecture 5.5} Let $p>3$ be a prime. Then
$$\align &\sum_{n=0}^{p-1}\f{(n+3)Q_n}{(-8)^n}
\e \cases 3p^2\mod{p^3}&\t{if $p\e 1\mod 3$,}
\\-15p^2\mod {p^3}&\t{if $p\e 2\mod 3$,}\endcases
\\& \sum_{n=0}^{p-1}\f{(n-2)Q_n}{(-9)^n}\e
 \cases -2p^2\mod{p^3}&\t{if $p\e 1\mod 3$,}
\\14p^2\mod {p^3}&\t{if $p\e 2\mod 3$.}\endcases\endalign$$
\endpro
\pro{Conjecture 5.6} Let $p$ be a prime with $p>3$. Then
$$\sum_{n=0}^{p-1}\f{Q_n}{(-6)^n}\e \cases
2x-\f p{2x}\mod {p^2}&\t{if $3\mid p-1$, $p=x^2+3y^2$ and $3\mid
x-1$},
\\-\f p{\b{(p-1)/2}{(p-5)/6}}\mod {p^2}&\t{if $p\e 2\mod 3$}
\endcases$$
and
$$\sum_{n=0}^{p-1}\f{Q_n}{(-12)^n}\e \cases
2x-\f p{2x}\mod {p^2}&\t{if $3\mid p-1$, $p=x^2+3y^2$ and $3\mid
x-1$},
\\\f p{2\b{(p-1)/2}{(p-5)/6}}\mod {p^2}&\t{if $p\e 2\mod 3$.}
\endcases$$
\endpro

\pro{Conjecture 5.7} Let $p$ be a prime with $p>3$. Then
$$\sum_{n=0}^{p-1}\b{2n}n\f{Q_n}{18^n}\e \cases
4x^2-2p-\f {p^2}{4x^2}\mod {p^3}&\t{if $p\e 1\mod 3$ and so
$p=x^2+3y^2$},
\\\f {p^2}{\b{(p-1)/2}{(p-5)/6}^2}\mod {p^3}&\t{if $p\e 2\mod 3$}
\endcases$$
and
$$\sum_{n=0}^{p-1}\b{2n}n\f{Q_n}{(-36)^n}\e \cases
4x^2-2p-\f {p^2}{4x^2}\mod {p^3}&\t{if $p\e 1\mod 3$ and so
$p=x^2+3y^2$,}
\\-\f {p^2}{2\b{(p-1)/2}{(p-5)/6}^2}\mod {p^3}&\t{if $p\e 2\mod 3$.}
\endcases$$
\endpro

\pro{Conjecture 5.8} Let $p>5$ be a prime. Then
$$\sum_{n=0}^{p-1}\b{2n}n\f{Q_n}{(-27)^n} \e
\sum_{n=0}^{p-1}\b{2n}n\f{Q_n}{(-81)^n}\e
\sum_{n=0}^{p-1}\b{2n}n\f{a_n}{(-45)^n} \e
\sum_{n=0}^{p-1}\b{2n}n\f{a_n}{9^n} \mod {p^2}.$$ If $p\e
1,17,19,23\mod{30}$, then
$$\align &\Ls {-3}p\sum_{n=0}^{p-1}\b{2n}n\f{a_n}{9^n}
\\& \e\cases 4x^2-2p-\f{p^2}{4x^2}\mod {p^3}&\t{if $p\e
1,19\mod {30}$ and so $p=x^2+15y^2$,}
\\2p-12x^2+\f{p^2}{12x^2}\mod {p^3}&\t{if $p\e
17,23\mod {30}$ and so $p=3x^2+5y^2$.}
\endcases\endalign$$
If $p\e 7,11,13,29\mod{30}$, then
$$ \Ls{-3}p\sum_{n=0}^{p-1}\b{2n}n\f{a_n}{9^n}
 \e\cases \f {31}{16}p^2\cdot 5^{[p/3]}\b{[p/3]}{[p/15]}^{-2}\mod
{p^3}&\t{if $p\e 7\mod{30}$,}
\\\f{31}4p^2\cdot 5^{[p/3]}\b{[p/3]}{[p/15]}^{-2}\mod {p^3}&\t{if $p\e
11\mod{30}$,}
\\ \f {31}{256}p^2\cdot
5^{[p/3]}\b{[p/3]}{[p/15]}^{-2}\mod {p^3}&\t{if $p\e 13\mod{30}$,}
\\\f {31}{64}p^2\cdot 5^{[p/3]}\b{[p/3]}{[p/15]}^{-2}\mod {p^3}&\t{if $p\e
29\mod{30}$.}
\endcases$$
\endpro

\pro{Conjecture 5.9} Let $p>3$ be a prime. Then
$$(-1)^{\f{p-1}2}\sum_{n=0}^{p-1}\b{2n}n\f{a_n}{54^n}\e \cases
4x^2-2p\mod {p^2}&\t{if $p=x^2+3y^2\e 1\mod 3$},
\\0\mod {p^2}&\t{if $p\e 2\mod 3$.}
\endcases$$
\endpro
\pro{Conjecture 5.10} Let $p$ be an odd prime. Then
$$\align &\sum_{n=0}^{p-1}\b{2n}n\f{Q_n}{(-32)^n}
\e\sum_{n=0}^{p-1}\b{2n}n\f{Q_n}{64^n}
 \\&\e\cases
4x^2-2p\mod {p^2}&\t{if $p=x^2+2y^2\e 1,3\mod 8$},
\\0\mod {p^2}&\t{if $p\e 5,7\mod 8$.}
\endcases\endalign$$
\endpro
\pro{Conjecture 5.11} Let $p>5$ be a prime. Then
$$\align&(-1)^{\f{p-1}2}\sum_{n=0}^{p-1}\b{2n}n\f{a_n}{20^n}
\e (-1)^{\f{p-1}2}\sum_{n=0}^{p-1}\b{2n}n\f{Q_n}{(-16)^n}
\\&\e \cases 4x^2-2p\mod{p^2}&\t{if $p\e 1,9\mod{20}$ and so
$p=x^2+5y^2$,}\\2x^2-2p\mod{p^2}&\t{if $p\e 3,7\mod{20}$ and so
$2p=x^2+5y^2$,}
\\0\mod{p^2}&\t{if $p\e 11,13,17,19\mod{20}$.}\endcases\endalign$$\endpro

\pro{Conjecture 5.12} Let $p>3$ be a prime. Then
$$\sum_{n=0}^{p-1}\b{2n}n\f{Q_n}{(-48)^n}\e
\cases 4x^2-2p\mod{p^2}&\t{if $12\mid p-1$ and so $p=x^2+9y^2$,}
\\2p-2x^2\mod {p^2}&\t{if $12\mid p-5$ and so $2p=x^2+9y^2$,}
\\0\mod {p^2}&\t{if $p\e 3\mod 4$.}
\endcases$$
\endpro
\pro{Conjecture 5.13} Let $p$ be a prime with $p>3$. If
$m\in\{-7,-25,-169,$ \newline $-1519,-70225,20,56,650,2450\}$ and
$p\nmid m(m-2)$, then
$$\sum_{n=0}^{p-1}\b{2n}n\f{Q_n}{(16(m-2))^n}\e
\Ls {m(m-2)}p\sum_{n=0}^{p-1} \b{2n}n\f{f_n}{(16m)^n}\mod {p^2}.$$
\endpro

\pro{Conjecture 5.14} Let $p$ be a prime with $p>3$. If
$m\in\{-112,-400,-2704,\newline -24304,-1123600\}$ and $p\nmid
m(m+4)$, then
$$\sum_{n=0}^{p-1}\b{2n}n\f{a_n}{(m+4)^n}\e
\Ls {m(m+4)}p\sum_{n=0}^{p-1} \b{2n}n\f{f_n}{m^n}\mod {p^2}.$$
\endpro


\begin{thebibliography}{CTYZ}
\bibitem [A] {} S. Ahlgren, {\it
  Gaussian hypergeometric series and combinatorial congruences,} in:
Symbolic computation, number theory, special functions, physics and
combinatorics (Gainesville, FI, 1999), pp.1-12, Dev. Math., Vol. 4,
Kluwer, Dordrecht, 2001, pp.1-12.

\bibitem [AO] {} S. Ahlgren and K. Ono, {\it A Gaussian hypergeometic series evaluation
and Ap\'ery number
 congruences}, J. Reine Angew. Math. {\bf 518}(2000), 187-212.

 \bibitem [AZ] {} G. Almkvist and W. Zudilin, {\it
 Differential equations, mirror maps
 and zeta values}, in: Mirror Symmetry
  AMS/IP Stud. Adv. Math.
  {\bf 38}(2006),
 International Press $\&$ Amer. Math. Soc., 481-515.

\bibitem [Ap] {} R. Ap\'ery, {\it Irrationalit\'e de $\zeta(2)$ et $\zeta(3)$},
 Ast\'erisque {\bf 61}(1979), 11-13.

 \bibitem [B1] {} F. Beukers, {\it Some congruences for the Ap\'ery numbers},
 J. Number Theory {\bf 21}(1985), 141-155.

\bibitem [B2] {} F. Beukers, {\it Another congruence for the Ap\'ery numbers},
 J. Number Theory {\bf 25}(1987), 201-210.

\bibitem [CCS]{} H.H. Chan, S. Cooper and F. Sica, {\it Congruences
satisfied by Ap\'ery-like numbers}, Int. J. Number Theory {\bf
6}(2010), 89-97.

\bibitem [CTYZ]{} H.H. Chan, Y. Tanigawa, Y. Yang and W. Zudilin,
 {\it New analogues of Clausen's identities arising from the theory
 of modular forms}, Adv. Math. {\bf 228}(2011), 1294-1314.

 \bibitem [CZ] {} H.H. Chan and W. Zudilin, {\it New representations
 for Ap\'ery-like sequences}, Mathematika {\bf 56}(2010), 107-117.



\bibitem [GMY] {} A. Gomez, D. McCarthy, D. Young, {\it Ap\'ery-like numbers and
families of newforms with complex multiplication}, Res. Number
Theory 2019, 5:5, pp.1-12.

\bibitem [G] {} H.W. Gould, {\it Combinatorial
Identities, A Standardized Set of Tables Listing 500 Binomial
Coefficient Summations}, West Virginia University, Morgantown, WV,
1972.

\bibitem [JV] {} F. Jarvis and H.A. Verrill,
{\it Supercongruences for the Catalan-Larcombe-French numbers},
Ramanujan J. {\bf 22}(2010), 171-186.

\bibitem [JS] {} X.J. Ji and Z.H. Sun, {\it Congruences for
Catalan-Larcombe-French numbers}, Publ. Math. Debrecen {\bf
90}(2017), 387-406.

\bibitem [LN] {} J.C. Liu and H.X. Ni, {\it Supercongruences for
Almkvist-Zudilin sequences}, arXiv:2004.0765.

\bibitem [M] {} E. Mortenson, {\it Supercongruences for truncated $\
_{n+1}F_n$ hypergeometric series with applications to certain weight
three newforms}, Proc. Amer. Math. Soc.{\bf 133}(2005), 321-330.

\bibitem [St] {} V.Strehl, {\it Binomial identities
-combinatorial and algorithmic aspects}, Discrete Math. {\bf
136}(1994), 309-346.


\bibitem [SB] {} J. Stienstra and F. Beukers, {\it On the
Picard-Fuchs equation and the formed Brauer group of certain
elliptic $K3$-surfaces}, Math. Ann. {\bf 271}(1985), 269-304.

\bibitem [S1] {} Z.H. Sun, {\it Congruences concerning Bernoulli
numbers and Bernoulli polynomials}, Discrete Appl. Math. {\bf
105}(2000), 193-223.

\bibitem [S2] {} Z.H. Sun, {\it Congruences involving Bernoulli and
Euler numbers}, J. Number Theory {\bf 128}(2008), 280-312.

\bibitem [S3] {}
Z. H. Sun, {\it Congruences concerning Legendre polynomials}, Proc.
Amer. Math. Soc. {\bf 139}(2011), 1915-1929.

\bibitem [S4] {} Z.H. Sun, {\it Identities and congruences for a new sequence}, Int. J.
Number Theory {\bf 8}(2012), 207-225.


\bibitem [S5] {}
Z. H. Sun, {\it Congruences concerning Legendre polynomials II}, J.
Number Theory {\bf 133}(2013), 1950-1976.

\bibitem [S6] {}
Z. H. Sun, {\it Congruences involving $\b{2k}k^2\b{3k}k$}, J. Number
Theory {\bf 133}(2013), 1572-1595.

\bibitem [S7] {}
Z. H. Sun, {\it Legendre polynomials and supercongruences}, Acta
Arith. {\bf 159}(2013), 169-200.


\bibitem [S8] {}
Z. H. Sun, {\it Generalized Legendre polynomials and related
supercongruences}, J. Number Theory {\bf 143}(2014), 293-319.


\bibitem [S9] {} Z. H. Sun, {\it Congruences for Domb and
Almkvist-Zudilin numbers}, Integral Transforms Spec. Funct. {\bf
26}(2015), 642-659.

\bibitem [S10] {} Z. H. Sun, {\it
Super congruences concerning Bernoulli polynomials}, Int. J. Number
Theory {\bf 11}(2015), 2393-2404.

\bibitem [S11] {} Z. H. Sun, {\it Note on super congruences modulo
$p^2$}, arXiv:1503.0341.

\bibitem [S12] {} Z. H. Sun, {\it Supercongruences involving
Bernoulli polynomials}, Int. J. Number Theory {\bf 12}(2016),
1259-1271.

\bibitem [S13] {} Z. H. Sun, {\it Supercongruences involving
Euler polynomials}, Proc. Amer. Math. Soc. {\bf 144}(2016),
3295-3308.

\bibitem [S14] {} Z.H. Sun, {\it Some further properties of even and
odd sequences}, Int. J. Number Theory {\bf 13}(2017), 1419-1442.

\bibitem [S15] {} Z.H. Sun, {\it Identities and congruences for
Catalan-Larcombe-French numbers}, Int. J. Number Theory {\bf
13}(2017), 835-851.

\bibitem [S16] {} Z.H. Sun, {\it Congruences for sums involving
Franel numbers}, Int. J. Number Theory {\bf 14}(2018), 123-142.

\bibitem [S17] {} Z.H. Sun, {\it Super congruences for two
Ap\'ery-like sequences}, J. Difference Equ. Appl. {\bf 24}(2018),
1685-1713.

\bibitem [S18] {} Z.H. Sun, {\it Congruences involving binomial
coefficients and Ap\'ery-like numbers}, Publ. Math. Debrecen {\bf
96}(2020), 315-346.

\bibitem [S19] {} Z.H. Sun, {\it Super congruences concerning binomial
coefficients and Ap¡äery-like numbers }, arXiv:2002.12072.

\bibitem [S20] {} Z.H. Sun, {\it
New congruences involving Ap\'ery-like numbers}, arXiv:2004.07172.


\bibitem [Su1] {} Z.W. Sun, {\it Super congruences and Euler
numbers}, Sci. China Math. {\bf 54}(2011), 2509-2535.

\bibitem [Su2] {} Z.W. Sun, {\it
On sums of Apery polynomials and related congruences}, J. Number
Theory {\bf 132}(2012), 2673-2699.

\bibitem [Su3] {} Z.W. Sun, {\it On sums involving products of three binomial
coefficients}, Acta Arith. {\bf 156}(2012), 123-141.

\bibitem  [Su4] {} Z.W. Sun, {\it Conjectures and results on $x^2$
mod $p^2$ with $4p=x^2+dy^2$}, in: Number Theory and Related Areas
(eds., Y. Ouyang, C. Xing, F. Xu and P. Zhang), Higher Education
Press $\&$ International Press, Beijing and Boston, 2013,
pp.149-197.

\bibitem [Su5] {} Z.W. Sun, {\it Some new series for $1/\pi$
 and related congruences}, Nanjing Univ. J.
Math. Biquarterly {\bf 31}(2014), 150-164.

\bibitem [Su6] {} Z.W. Sun, {\it New series for some special values
of L-functions}, Nanjing Univ. J. Math. Biquarterly {\bf 32}(2015),
 189-218.

\bibitem [Su7] {} Z.W. Sun, {\it
Supercongruences involving dual sequences}, Finite Fields Appl. {\bf
46}(2017), 179-216.


\bibitem [T1] {} R. Tauraso, {\it Supercongruences for a truncated hypergeometric
series}, Integers {\bf 12}(2012), A45,12pp.

\bibitem [T2]{} R. Tauraso, {\it Supercongruences related to $\ _3F_2(1)$
 involving harmonic numbers}, Int. J. Number Theory {\bf 14}(2018),
 1093-1109.

 \bibitem [W] {} C. Wang, {\it Supercongruences and hypergeometric
 transformations}, arXiv:2003.09888.

\bibitem [Z]{} D. Zagier,
{\it Integral solutions of Ap\'ery-like recurrence equations}, in:
Groups and Symmetries: From Neolithic Scots to John McKay, J. Harnad
and P. Winternitz, eds., CRM Proceedings and Lecture Notes, Vol. 47
(American  Mathematical Society, Providence, RI, 2009), pp.349-366.

\bibitem [Zu] {} W. Zudilin, {\it Ramanujan-type formulae for $1/\pi$:
 A second wind?}, in: Modular Forms and String Duality
 (Banff, June 3-8, 2006), N. Yui, H. Verrill, and C.F. Doran (eds.),
 Fields Inst. Commun. Ser.
 {\bf 54}(2008), Amer. Math. Soc. $\&$ Fields Inst., 179-188.




\end{thebibliography}
\end{document}